\documentclass{article}

\usepackage[centertags]{amsmath}

\usepackage{amssymb}
\usepackage{latexsym}
\usepackage{amscd}
\usepackage{theorem}
\usepackage{epic}
\usepackage{eepic}
\usepackage{epsfig}

{\theorembodyfont{\slshape}

        \newtheorem{thm}{Theorem}[section]
        \newtheorem{cor}[thm]{Corollary}
        \newtheorem{lem}[thm]{Lemma}
        \newtheorem{prop}[thm]{Proposition}
}

{\theorembodyfont{\rmfamily}

        \newtheorem{defn}[thm]{Definition}
        \newtheorem{prob}[thm]{Problem}
        \newtheorem{rem}[thm]{Remark}

        \newtheorem{claim}[thm]{Claim}
        
        \newtheorem{question}[thm]{Question}

}

\makeatletter
\renewcommand{\subsection}{\@startsection{subsection}{3}%
        {\z@}{-3.25ex plus -1ex minus-.2ex}{-1em}{\bf}}
\makeatother

\newcommand{\proof}{\noindent {\bf Proof:\ }}
\newcommand{\Endproof}{\hspace*{\fill} $\Box$ \vspace{1ex} \noindent }

\newcommand{\ZZ}{\mathbb{Z}}

\newcommand{\PP}{\mathbb{P}}

\newcommand{\FF}{\mathbb{F}}

\renewcommand{\AA}{\mathbb{A}}

\newcommand{\Hc}{{\mathcal H}}
\newcommand{\C}{{\mathcal C}}

\newcommand{\T}{\mathcal T}

\newcommand{\Sb}{\bar{S}}
\newcommand{\Tb}{\bar{T}}
\newcommand{\Ub}{\bar{U}}

\newcommand{\Spec}{\mathop{\rm Spec}} 

\newcommand{\Spf}{\mathop{\rm Spf}}

\newcommand{\Aut}{\mathop{\rm Aut}\nolimits}

\newcommand{\Def}{\mathop{\rm Def}\,}

\newcommand{\ord}{{\rm ord}}
\newcommand{\Ord}{{\rm ord}}

\newcommand{\res}{\mathop{\rm res}\nolimits}
\newcommand{\Proj}{\mathop{\rm Proj}}

\newcommand{\inj}{\hookrightarrow}
\newcommand{\To}{\;\longrightarrow\;}
\newcommand{\iso}{\stackrel{\sim}{\to}}

\newcommand{\NN}{\mathbb{N}}

\title{The local lifting problem for dihedral groups}

\date{}

\author{ Irene I.\ Bouw \and Stefan Wewers}

\begin{document}

\maketitle


\begin{abstract}
Let $G=D_p$ be the dihedral group of order $2p$, where $p$ is an odd
prime. Let $k$ an algebraically closed field of
characteristic $p$. We show that any action of $G$ on the ring
$k[[y]]$ can be lifted to an action on $R[[y]]$, where $R$ is some
complete discrete valuation ring with residue field $k$ and fraction
field of characteristic $0$. 


\end{abstract}

\section{Introduction}

\subsection{} 

Let $k$ be an algebraically closed field of characteristic $p$ and $G$
a finite group. A {\sl local $G$-action} is a faithful $k$-linear
action $\phi:G\inj\Aut_k(k[[y]])$.

\begin{prob}[The local lifting problem]\label{LLP}
  Does there exist a lift of $\phi$ to an $R$-linear action
  $\phi_R:G\inj\Aut_R(R[[y]])$, where $R$ is some complete
  discrete valuation ring with residue field $k$ and fraction field of
  characteristic zero?
\end{prob}
If a lift $\phi_R$ as above exists then we say that $\phi$ {\sl lifts
  to characteristic $0$}.

The main motivation for studying Problem \ref{LLP} comes from the
global lifting problem: given a smooth projective curve $Y$ over $k$
on which the group $G$ acts, can we lift the curve $Y$ to
characteristic $0$, together with the $G$-action? It is well known
that it suffices to construct such a lift locally, i.e.\ to lift the
completion of the curve $Y$ at closed points together with the action
of the stabilizer of that point. The global lifting problem is hence
reduced to the local lifting problem. See \cite{GM98} or
\cite{BM00}.  In this paper we solve the local lifting problem
for the dihedral group of order $2p$.

\begin{thm} \label{thm0}
  Suppose that $p$ is odd and that $G$ is the dihedral group of order
  $2p$. Then every local $G$-action $\phi:G\inj\Aut_k(k[[y]])$ lifts
  to characteristic $0$.
\end{thm}

Our approach to the local lifting problem is based on a generalization
of the methods of Green--Matignon \cite{GM99} and Henrio
\cite{Henrio02} (which treat the case $G=\ZZ/p$). We roughly do the
following. Let $G=\ZZ/p\rtimes\ZZ/m$ be a semidirect product of a
cyclic group of order $m$ by a cyclic group of order $p$ (with
$(p,m)=1$). Suppose that $G$ acts on the open rigid disc over a
$p$-adic field $K$. To this action we associate a certain object,
called a {\sl Hurwitz tree}. This object represents, in some sense, the
reduction of the $G$-action on the disc to characteristic $p$.  Our
first main result is that one can reverse this construction. As a
consequence, the local lifting problem for the group
$G=\ZZ/p\rtimes\ZZ/m$ can be reduced to the construction of certain
Hurwitz trees. Our second main result is that this can always be done
if $G=D_p$. This involves the construction of certain differential
forms on the projective line with very specific properties.

\subsection{}

Let $G$ be a finite group. If there exists a local $G$-action
$\phi:G\inj\Aut_k(k[[y]])$, then $G$ is an extension of a cyclic group
$C$ of order prime to $p$ by a $p$-group $P$. Furthermore, if $P=1$
then $\phi$ always lifts to characteristic $0$.

In this paper we  always assume that the subgroup $P$ has order
$p$. Then there exists a character $ \chi:C\to \FF_p^\times$ such that
$G=P\rtimes_\chi C$ (by this we mean that
$\tau\sigma\tau^{-1}=\sigma^{\chi(\tau)}$ for $\sigma\in P$ and
$\tau\in C$). Write $m$ for the order of $C$. 

Let $\phi:G\inj\Aut_k(k[[y]])$ be a local $G$-action. Let $\sigma$ be
a generator of the cyclic group $P$. The {\sl conductor} of $\phi$ is
the positive integer
\[
     h  :=  \ord_y (\frac{\sigma(y)}{y}-1).
\] 
There are certain obvious necessary conditions on the character $\chi$
and the conductor $h$ for a lift to exist. In the case where $G=P$ has
order $p$, this was first noted by Oort, \cite[\S I.1]{Oort85}.  In
\cite{Bertin}, Bertin found a systematic way to produce necessary
conditions for the liftability of $\phi$, for general groups $G$. This
will be further discussed in \cite{CGH}. In our case (where $P$ has
order $p$), we can summarize these conditions as follows (see
Corollary \ref{Htreecor1} below).

\begin{prop}\label{NC}
  Suppose that $\phi$ lifts. Then the following holds.
  \begin{enumerate}
  \item
    The character $\chi$ is either trivial or injective. 
  \item
    If $\chi$ is injective, $m|h+1$. 
  \end{enumerate}
\end{prop}

Note that for the group $G=D_p$, $\chi$ is injective and $m=2$.
Moreover, the conductor of any local $D_p$-action $\phi$ is odd, by
the Hasse--Arf theorem. Hence the necessary conditions in Proposition
\ref{NC} are automatically verified.

Chinburg--Guralnick--Harbater \cite{CGH} call a group $G=P\rtimes_\chi
C$ a {\sl local Oort group at $p$} if every local $G$-action
$\phi:G\inj\Aut_k(k[[y]])$ lifts to characteristic zero. This
terminology is inspired by the conjecture of Oort \cite{Oort95} that
every cyclic cover of smooth projective curves is liftable to
characteristic zero. Green--Matignon \cite{GM98} have shown that a
cyclic group of order $pm$, with $(p,m)=1$, is a local Oort group at
$p$.  Proposition \ref{NC} shows that, for $P$ cyclic of order $p$,
the only other possible local Oort group at $p$ is the dihedral group
$D_p$ of order $2p$. Theorem \ref{thm0} shows that $D_p$ is indeed a
local Oort group. Hence we now have a complete list of all local Oort
groups with Sylow $p$-subgroup of order $p$.

\begin{question} \label{question1}
  Let $G=P\rtimes_\chi C$ be as above, with $P\cong\ZZ/p$, and let
  $\phi:G\inj\Aut_k(k[[y]])$ be a local $G$-action. Suppose that $\phi$
  satisfies the necessary conditions given by Proposition \ref{NC}. Is
  it then true that $\phi$ always lifts?
\end{question}

In this paper we did not attempt to answer Question \ref{question1}.
Nevertheless, the methods presented here reduce this question to
solving certain explicit equations over $\FF_p$.

\subsection{}

The local lifting problem (Problem \ref{LLP}) makes sense also for
groups $G$ which are the extension of a cyclic group of order prime to
$p$ by an arbitrary $p$-group. For $p$ odd and $G\cong(\ZZ/p)^2$,
Green and Matignon have constructed local $G$-actions $\phi$ which do
not lift. They also showed that every local $G$-action lifts if $G$ is
a cyclic group of order $p^nm$ with $n=1,2$ and $\gcd(p,m)=1$, see
\cite{GM98}.  Therefore, it seems reasonable to conjecture that for
$p$ odd the only local Oort groups at $p$ are the cyclic groups of
order $p^nm$ and the dihedral groups of order $2p^n$. 

For $p=2$ the situation is somewhat different. The groups $(\ZZ/2)^2$
(\cite[Th.\ 5.2.1]{Pagot}) and $A_4$ (I.I.\ Bouw, unpublished) are
local Oort groups at $2$.  It seems one needs a more careful analysis
of the situation before one can make a reasonable guess on what to
expect here.

\vspace{2ex} The paper is organized as follows. In \S \ref{bound} we
study $G$-action on the boundary of a rigid disc. This study is
extended in \S \ref{Htree} to $G$-actions on the whole disc. The main
result here is Theorem \ref{Htreethm} which allows the construction of
a $G$-action on the disc, starting from a Hurwitz tree.  In \S
\ref{largehsec} we prove our main theorem (Theorem \ref{thm0}) in the
case where the conductor $h$ is $>p$. In \S \ref{smallhsec} we deal
with the case $h<p$.

\vspace{2ex} We thank Leonardo Zapponi for interesting comments on the
previous version of this paper. Among other things, he told us how the
equations in \S \ref{trivsolsec} could be simplified.


\section{Group actions on the boundary of the disc} \label{bound}

In this section we study possible actions of the group
$G=\ZZ/p\rtimes\ZZ/m$ on the boundary of an open rigid disc. This is a
preparation for the next section, where we will study $G$-actions on
the whole disc.  The main result is Proposition \ref{boundprop3} which
shows that the action of $G$ on the boundary of a disc is determined,
up to conjugation, by three simple invariants. This will allow us in
\S \ref{Htree} to construct $G$-actions on the whole disc by
patching. In the case $G=\ZZ/p$ this is already contained in
\cite{Henrio02}.

\subsection{} \label{bound1}

Let $R$ be a complete discrete valuation ring with residue field $k$
and fraction field $K$. We assume that $k$ is algebraically closed of
characteristic $p>2$ and that $K$ has characteristic $0$. We fix a
uniformizing element $\pi$ of $R$.

Let 
\[
      Y_K  :=  \{\,y\,\mid\, |y|_K < 1\,\}.
\]
be the rigid open unit disc with parameter $y$.  We let
$S=R[[y]]\{y^{-1}\}$ denote the ring of formal Laurent series
$f=\sum_i f_i\,y^i$ with $f_i\in R$ for $i\in\ZZ$ and $\lim_{i\to
-\infty} f_i=0$. The ring $S$ is a complete discrete valuation ring
with residue field $\Sb=k((y))$, uniformizing element $\pi$ and
fraction field $S\otimes_R K$. We can think of elements of $S$ as
bounded functions on the boundary of $Y_K$.

We denote by $v:(S\otimes
K)^\times\to\ZZ$ the  exponential valuation normalized by $v(\pi)=1$ and
by $\ord_y:\Sb^\times\to\ZZ$ the exponential valuation normalized by 
$\ord_y(y)=1$. Furthermore, we define the function $\natural:(S\otimes
K)^\times\to\ZZ$ by the formula
\[
        \natural f  :=  \ord_y(\frac{f}{\pi^{v(f)}}\mod{\pi}).
\]

A {\sl parameter} on the boundary of the disc is a unit $z\in
S^\times$ with $\natural z=1$. In other words, $z=\sum_i a_iy^i$ with
$a_1\in R^\times$ and $a_i\in (\pi)$ for $i<1$. If $z$ is a parameter,
then $S=R[[z]]\{z^{-1}\}$, i.e.\ every element of $S$ can be written
as a formal Laurent series $f=\sum_i b_i\,z^i$ with $b_i\in R$ for
$i\in\ZZ$ and $\lim_{i\to -\infty} b_i=0$.

An {\sl automorphism} of the boundary of the disc is an automorphism
of $\Spec S$ which induces a continuous automorphism of the residue
field $\Sb=k((y))$ of $S$. We write $\Aut_R(\Spec S)$ for the group of
such automorphisms. If $z$ is any parameter,  there exists a
unique automorphism $\tau$ with $\tau^*y=z$.

\subsection{} \label{bound2}

Let $P$ be a cyclic group of order $p$, $C$ a cyclic group of order
$m$, with $(m,p)=1$, and $\chi:C\to\FF_p^\times$ a character. We set
$G:=P\rtimes_\chi C$.

A $G$-action on the boundary of the disc is an injective group
homomorphism $\psi:G\inj\Aut_R(\Spec S)$. In this subsection we define
three invariants of such an action: the { conductor}, the { different}
and the {tame inertia character}. In the next subsection we will show
that $\psi$ is uniquely determined, up to conjugation in $\Aut_R(\Spec
S)$, by these invariants. This is a generalization of 
\cite[Corollaire 1.8]{Henrio02}, which treats the case $G=P$.

Fix a $G$-action $\psi:G\inj\Aut_R(\Spec S)$, and consider the group
$G$ as a subgroup of $\Aut_R(\Spec S)$ by means of $\psi$. Let
$T:=S^P$ be the ring of invariants under the action of $P\subset G$.
Clearly, $T$ is again a complete discrete valuation ring with
parameter $\pi$, and $S$ is finite and flat over $T$ of degree $p$. It
follows from classical valuation theory that the extension of residue
fields $\Sb/\Tb$ has degree $p$ and that $\Tb=k((z))$ is a field of
Laurent series over $k$. Let $\mathcal{D}_{S/T}$ denote the different
of the extension $S/T$.

\begin{defn} \label{bounddef}
  The {\sl different} of the $G$-action $\psi$ on the boundary of the
  disc is the unique rational number $\delta$ with
  \[
     \mathcal{D}_{S/T} = p^\delta\cdot S.
  \]
  The {\sl conductor} of $\psi$ is the integer
  \[
       h  :=  \natural(\frac{\sigma^*y}{y}-1).
  \]   
  The {\sl tame inertia character} of $\psi$ is the character
  $\lambda:C\to k^\times$ such that
  \[
       \tau^*y  \equiv  \lambda(\tau)\cdot y \pmod{y^2}
  \]
  in $\Sb$. 
\end{defn}

The conductor and the different are already defined by
Henrio \cite{Henrio02}. He shows that  $0\leq\delta\leq 1$. Moreover,
we have the following classification result. 

\bigskip
{\bf The \'etale case}:
    Suppose that $\delta=0$. Then the extension $S/T$ is unramified
    and the extension of residue fields $\Sb/\Tb$ is Galois,
    with Galois group $P$. The conductor $h$ of $\psi$ is equal to
    the usual conductor of the Galois extension $\Sb/\Tb$. In
    particular, $h$ is prime to $p$ and positive.

\bigskip
{\bf The additive case}: 
 Suppose that $0<\delta<1$. Then the group $P$, considered
    as constant group scheme over $K$, extends to a finite and flat
    group scheme $P_R$ over $R$, with special fiber isomorphic to
    $\alpha_p$. Moreover, the action of $P$ on $\Spec(S\otimes_RK)$
    extends uniquely to an action of $P_R$ on $\Spec S$ which makes
    the map $\Spec S\to\Spec T$ a $P_R$-torsor. There exists an exact
    differential $\omega\in\Omega_{\Tb/k}^{\rm cont}$ which classifies
    the induced $\alpha_p$-torsor $\Spec\Sb\to\Spec\Tb$
    (\cite[Proposition III.4.14]{Milne}). We have
    \[
         h = -\ord_z\omega - 1.
    \]
    Moreover, $h$ is prime to $p$ and different from $0$.
 
\bigskip
{\bf The multiplicative case}:
 Suppose that $\delta=1$. Then the group $P$, considered
    as constant group scheme over $K$, extends to a the finite and
    flat group scheme $P_R\cong\mu_{p,R}$ over $R$. Moreover, the
    action of $P$ on $\Spec(S\otimes_RK)$ extends uniquely to an
    action of $P_R$ on $\Spec S$ which makes the map $\Spec S\to\Spec
    T$ a $P_R$-torsor.  There exists a logarithmic differential
    $\omega\in\Omega_{\Tb/k}^{\rm cont}$ which classifies the induced
    $\mu_p$-torsor $\Spec\Sb\to\Spec\Tb$ (\cite[Proposition
    III.4.14]{Milne}). Again we have
    \[
         h = -\ord_z\omega - 1.
    \]
    Moreover, $h$ is prime to $p$ (except if $h=0$) and less then or
    equal to $0$.

See \cite[\S 1]{Henrio02} for details. Note that the invariants $h$
and $\delta$ only depend on the action of $P$. However, the fact that
the action of $P$ extends to an action of $G=P\rtimes_\chi C$ puts
some extra conditions on the conductor $h$.

\begin{prop} \label{boundprop2}
  The tame inertia character $\lambda$ has order $m=|C|$, and we have
  \[
       \chi = \lambda^{-h}.
  \]
  Therefore, the order of $h$ in $\ZZ/m$ is equal to the order of the
  character $\chi$. In particular:
  \begin{enumerate}
  \item
    If $\chi$ is injective then $h$ and $m$ are relatively
    prime. Moreover, $\lambda$ takes values in $\FF_p^\times$ and is
    uniquely determined by $\chi$ and $h$.
  \item
    If $\chi$ is trivial then $m|h$. 
  \end{enumerate}
\end{prop}

\proof
Let $U:=S^G$. As before, we conclude that $S/U$ is a Galois extension
of discrete valuation rings, with Galois group $G$ and residue field
extension $\Sb/\Ub$ of degree $|G|=pm$. Since $m$ is prime to $p$,
$\Tb/\Ub$ is actually a Galois extension with Galois group $C$. 

In the \'etale case $\Sb/\Ub$ is a Galois
extension with Galois group $G$. The filtration of higher ramification
groups in the upper numbering has a unique jump at $\sigma=h/m$.  In
this case, the statements of Proposition \ref{boundprop2} are well
known. For instance, (i) and (ii) follow from the Hasse--Arf theorem.
We recall the argument, because this allows us to introduce some
notation which we  need later on. By Artin--Schreier
theory, there exists a generator $w$ of the extension $\Sb/\Tb$ such
that $\sigma^*w=w+1$ and
\begin{equation} \label{ASeq1}
     w^p-w = z^{-h}+c_1z^{-h+1}+\ldots \in\Tb=k((z)).
\end{equation}
Also, it is a standard fact on tame ramification that we may choose
the parameter $z$ in such a way that $\tau^*z=\lambda'(\tau)\cdot z$
for all $\tau\in C$ and for some character $\lambda':C\to k^\times$ of
order $m$. It follows from \eqref{ASeq1} that $\ord_y(w)=-h$, where
$y$ is a parameter of $S$.
Therefore, we may assume (by changing the parameter $y$ of $\Sb$) that
$w=y^{-h}$. A straightforward computation, using
$\tau\sigma\tau^{-1}=\sigma^{\chi(\tau)}$ shows that
$\chi(\tau)=\lambda'(\tau)^{-h}$ and $\tau^*y\equiv\lambda'(\tau)\cdot
y\pmod{y^2}$. Hence $\lambda=\lambda'$ and $\lambda^{-h}=\chi$. This
proves the proposition in the \'etale case.

In the additive and multiplicative case the extension $\Sb/\Tb$ is
inseparable of degree $p$ and an $\alpha_p$- or $\mu_p$-torsor. The
torsor structure is encoded by a differential form
$\omega\in\Omega_{\Tb/k}$ with $\ord_z\omega=-h-1$. If $w^p=u$ is an
equation for this torsor, then $\omega={\rm d}u$ (in the case of an
$\alpha_p$-torsor) or $\omega={\rm d}u/u$ (in the case of a
$\mu_p$-torsor). See \cite{Henrio02} for details. It is easy to check
that the rule $\tau\sigma\tau^{-1}=\sigma^{\chi(\tau)}$ implies
$\tau^*\omega=\chi(\tau)\cdot\omega$. As in the \'etale case, we
choose a parameter $z$ of $\Tb$ with $\tau^*z=\lambda'(\tau)\cdot z$
for a character $\lambda'$ of order $m$. Then
$\omega=(az^{-h-1}+\ldots)\,{\rm d}z$, and hence
$\chi=(\lambda')^{-h}$. If $h\not=0$ then we can choose the
torsor-equation $w^p=u$ such that $u=1+z^{-h}+\ldots$. It follows that
$y:=(w-1)^{-1/h}$ is a parameter for $\Sb$ with
$\tau^*y\equiv\lambda'(\tau)\cdot y\pmod{y^2}$. Hence
$\lambda=\lambda'$ and $\lambda^{-h}=\chi$. On the other hand, if
$h=0$ then we see that $\chi$ must be trivial, and so the claim of the
proposition is true as well. This finishes the proof of the
proposition.  \Endproof

\subsection{} \label{bound3}

Let $\psi_i:G\inj\Aut_R(S)$, $i=1,2$, be two $G$-actions on the
boundary of the disc. We say that $\psi_1$ and $\psi_2$ are {\sl
  conjugate} if there exists an automorphism $\eta\in\Aut_R(\Spec S)$
such that $\psi_2(g)\circ\eta=\eta\circ\psi_1(g)$ holds for all $g\in
G$.

\begin{prop} \label{boundprop3}
  Suppose that the two $G$-actions $\psi_1$ and $\psi_2$ have the same
  conductor $h$, the same different $\delta$ and the same tame inertia
  character $\lambda$. Suppose, moreover, that $h\not=0$. Then
  $\psi_1$ and $\psi_2$ are conjugate.
\end{prop}

\proof It follows from \cite[Corollaire 1.8]{Henrio02} that the
restrictions of $\psi_1$ and $\psi_2$ to the subgroup $P\subset G$ are
conjugate. Hence we may assume that $\psi_1|_P=\psi_2|_P$.

The $P$-action on $\Spec(S\otimes K)$
extends to a free action of a certain $R$-group scheme $P_R$ on $S$,
where $P_R$ is either the constant group scheme $P$ (\'etale case) or a
connected group scheme with special fiber $\alpha_p$ (additive case) or
$\mu_p$ (multiplicative case). For $n\geq 0$, set
$S_n:=S\otimes_R(R/\pi^{n+1})$ and $P_n:=P_R\otimes_R(R/\pi^{n+1})$.

\begin{claim} \label{boundprop3claim}
  There exists a family $(\eta_n)$, where $\eta_n$ is an automorphism
  of $\Spec S_n$ which commutes with the action of $P_n$ and such that
  $\psi_2(\tau)\circ\eta_n\equiv\eta_n\circ\psi_1(\tau)$ and
  $\eta_{n+1}\equiv\eta_n\pmod{\pi^{n+1}}$.
\end{claim}

The claim implies that $\eta:=\lim_n\eta_n$ is an automorphism of $S$
which commutes with the action of $P$ and such that
$\psi_2(\tau)\circ\eta=\eta\circ\psi_1(\tau)$ for all $\tau\in C$. Hence
the proposition is an immediate consequence of Claim
\ref{boundprop3claim}.

We prove Claim \ref{boundprop3claim} by induction on $n$, and we start
at $n=0$. We consider the \'etale case first. Let $w$ be an Artin--Schreier
generator of the extension $\Sb/\Tb$, as in the proof of Proposition
\ref{boundprop2}. Then $\sigma^*w=w+1$ and
$f:=w^p-w=z^{-h}+\ldots\in\Tb$. Replacing $w$ by $w+g$, for a
suitable $g\in k[[z]]\cdot z^{-h+1}$, we can achieve that $f$ has no
monomials whose exponents are divisible by $p$. Furthermore, choosing
the parameter $z$ of $\Tb$ appropriately, we may assume that
$\psi_1(\tau)^*z=\lambda(\tau)\cdot z$ for all $\tau\in C$. It follows
from Artin--Schreier theory that
\[
    \psi_1(\tau)^*w = \chi(\tau)(w+g_\tau),
\]
with $g_\tau\in k[[z]]\cdot z^{-h+1}$, and hence
\[
     \psi_1(\tau)(f)^* = \chi(\tau)(f+g_\tau^p-g_\tau).
\]
Using our assumption on $f$ and $\psi_1(\tau)^*z=\lambda(\tau)\cdot
z$, it is now easy to see that $g_\tau=0$. We conclude that, in terms
of the parameter $y=w^{-h}$ for $\Sb$, the reduction of $\psi_1$ to
$\Sb$ is given by 
\[
   \psi_1(\sigma)^*y = y\,(1+y^h)^{-1/h}, \qquad
   \psi_1(\tau)^*y = \lambda(\tau)\cdot y.
\]
This depends visibly only on $h$ and $\lambda$. Hence we can find
another parameter $y'$ for $\Sb$ such that the reduction of $\psi_2$
to $\Sb$ is given by the same formulas (but with $y$ replaced by
$y'$). Now $\eta_0^*y:=y'$ defines the desired automorphism
$\eta_0$. This proves Claim \ref{boundprop3claim} for $n=0$.

The proof for $n=0$ in the additive and the multiplicative case is
quite similar. We sketch it in the multiplicative case. Then $\Sb\to\Tb$ is a
$\mu_p$-torsor, given by a Kummer equation $w^p=f$, with
$f\in\Tb^\times$. We choose $z$ such that
$\psi_1(\tau)^*z=\lambda(\tau)\cdot z$. Since $h<0$, we can choose
$f=1+z^{-h}+\ldots$, and may assume that $f$ has no monomial whose
exponent is divisible by $p$. Then $\psi_1(\tau)^*w=w^\nu g_\tau$,
with $\nu\equiv\chi(\tau)\pmod{p}$ and $g_\tau\in
k[[z]]^\times$. Therefore, $\psi_1(\tau)^*f=f^\nu g^p_\tau$. An
argument as before shows that $g_\tau=1$. The conclusion is as in the
\'etale case. This completes the proof of Claim \ref{boundprop3claim} for
$n=0$.

Let us assume that we have constructed $\eta_n$ for some $n\geq 0$.
We claim that we can lift $\eta_n$ to an automorphism $\eta'$ of $S$
which commutes with the action of $P_R$. For instance, in the
multiplicative case we have $\eta_n(y)=y v_n$, where $v_n\in
T_{n}=S_{n}^{P_n}$ and $y^p=f$ is a Kummer equation for the
$\mu_p$-torsor $\Spec S\to\Spec T$. Then we can simply lift $v_n$ to
an element $v\in T$ and set $\eta'(y):=yv$.  The construction of the
lift $\eta'$ in the \'etale and the additive case is similar.

Let us choose a lift $\eta'$ as above. Conjugating $\psi_2$ by
$\eta'$, we may assume that $\psi_1\equiv\psi_2\pmod{\pi^{n+1}}$. Then
for $\tau\in C$ we define
$\alpha_\tau:=\psi_2(\tau)\circ\psi_1(\tau)^{-1}\mod{\pi^{n+2}}$.
This is an automorphism of $\Spec S_{n+1}$ which is the identity on
$\Spec S_n$. Therefore, for all $w\in S_{n+1}$ we have
\[
     \alpha_\tau^*w = w + \pi^{n+1}\cdot D_\tau(\bar{w}),
\]
where $D_\tau$ is a (continuous) $k$-derivation of $\Sb$. It is clear
that $\alpha_\tau$ commutes with the action of $P_{n+1}$. Therefore
the derivation $D_\tau$ is invariant under the action of $P_0$. One
also checks that the map
\[
   C  \To  {\rm Der}_{\Sb/k}^{P_0}, \qquad
   \tau  \mapsto  D_\tau
\]
is a cocycle (note that ${\rm Der}_{\Sb/k}^{P_0}$ is a $C$-module in a
natural way). Since ${\rm Der}_{\Sb/k}^{P_0}$ is a $k$-vector space and
the order of $C$ is prime to $p$, this cocycle is a coboundary. This
means that there exists a $P_0$-invariant derivation $D$ such that
$D_\tau=\tau\circ D\circ\tau^{-1}-D$. Set
$\eta_{n+1}^*w:=w-\pi^{n+1}\cdot D(\bar{w})$. Then
$\psi_2(\tau)\circ\eta_{n+1}\equiv\eta_{n+1}\circ\psi_1(\tau)$ for all
$\tau$, as desired. By induction, we get a proof of Claim
\ref{boundprop3claim}, and hence of Proposition \ref{boundprop3}.
\Endproof

\begin{rem} \label{boundrem1}
  The proof of Proposition \ref{boundprop3} in the \'etale case for
  $n=0$ shows the following. Let $\phi_1,\phi_2:G\inj\Aut_k(k[[y]])$
  be two local $G$-actions with the same conductor and the same tame
  inertia character. Then $\phi_1$ and $\phi_2$ are conjugate (inside
  $\Aut_k(k[[y]])$).
\end{rem}

\section{Group actions on the disc} \label{Htree}

In this section we study actions of the group $G=\ZZ/p\rtimes\ZZ/m$ on
the open rigid disc. We generalize the results
of Green--Matignon \cite{GM99} and Henrio \cite{Henrio02}, which treat
the case $G=\ZZ/p$. 

Throughout, we suppose that $p$ is an odd prime and that $k$ is an
algebraically closed field of characteristic $p$.

\subsection{Hurwitz trees}\label{Htree1}

We start by giving a formal definition of a Hurwitz tree, following
Henrio \cite{Henrio02}. Our definition differs slightly from the
definition given in \cite{Henrio02}; for instance, we do not
distinguish between Hurwitz trees and realizations of Hurwitz trees.

\begin{defn} \label{dectreedef}
  A {\sl decorated tree} is given by the following data:
  \begin{itemize}
  \item
    a semistable curve $Z$ over $k$ of genus $0$,
  \item
    a family $(z_b)_{b\in B}$ of pairwise distinct smooth $k$-rational
    points of $Z$, indexed by a finite {\sl nonempty} set $B$,
  \item
    a distinguished smooth $k$-rational point $z_0\in Z$, distinct from
    any of the points $z_b$.
  \end{itemize}
  We require that $Z$ is {\sl stably marked} by the points
  $(z_b;z_0)$, see \cite{Knudsen83}. 

  The {\sl combinatorial tree} underlying a decorated tree $Z$ is the
  graph $T=(V,E)$, defined as follows. The vertex set $V$ of $T$ is
  the set of irreducible components of $Z$, together with
  a distinguished element $v_0$, called the {\sl root}. The edge set $E$
  is the set of singular points of $Z$, together with a distinguished
  element $e_0$. We write $Z_v$ for the component corresponding to a
  vertex $v\not=v_0$ and $z_e$ for the singular point corresponding to
  an edge $e\not=e_0$. An edge $e$ corresponding to a singular point
  $z_e$ is adjacent to the vertices corresponding to the two
  components which intersect in $z_e$. The edge $e_0$ is adjacent to
  the root $v_0$ and the vertex $v$ corresponding to the (unique)
  component $Z_v$ containing the distinguished point $z_0$.  For each
  edge $e\in E$, the {\sl source} (resp.\ the {\sl target}) of $e$ is
  the unique vertex $s(e)\in V$ (resp.\ $t(e)\in V$) adjacent to $e$
  which lies in the direction of the root (resp.\ in the direction
  away from the root).
\end{defn}

Note that, since $(Z,(z_b),z_0)$ is stably marked of genus $0$, the
components $Z_v$ have genus zero, too, and the graph $T$ is a tree.
Moreover, we have $|B|\geq 2$.

An {\sl admissible action} of a group $C$ on a decorated tree
$(Z,(z_b),z_0)$ is a $k$-linear action of $C$ on the curve $Z$ which
satisfies the following. Firstly, the action permutes the points $z_b$
and fixes the point $z_0$. Secondly, for any singular point $z_e$,
$e\not=e_0$, the stabilizer $C_e\subset C$ of $z_e$ is of order
prime-to-$p$, and the two characters describing the action of $C_e$ on
the tangent spaces of the two branches of $Z$ at the point $z_e$ are
inverse to each other. For every edge $e$ (including $e_0$) we shall
write $\lambda_e:C_e\to k^\times$ for the character which describes
the action of $C_e$ on the tangent space of the component $Z_{t(e)}$
at $z_e$.

For a vertex $v\in V$, we write $U_v\subset Z_v$ for the complement in
$Z_v$ of the set of singular and marked points.

\begin{defn} \label{Htreedef}
  Let $C$ be a cyclic group of order $m$ (with $(m,p)=1$) and
  $\chi:C\inj\FF_p^\times$ a character.  A {\sl Hurwitz tree} of type
  $(C,\chi)$ is given by the following data:
  \begin{itemize}
  \item
    A decorated tree $Z=(Z,(z_b),z_0)$ with underlying combinatorial
    tree $T=(V,E)$, together with an admissible action of $C$.
  \item
    For each $v\in V-\{v_0\}$, a differential form
    $\omega_v$ on $U_v\subset Z_v$ without zeroes or poles.
  \item
    For every $v\in V$, a rational
    number $0\leq\delta_v\leq 1$, called the {\sl different} of $v$.
  \item
    For every $e\in E$, a positive rational number $\epsilon_e$,
    called the {\sl thickness} of $e$.
  \end{itemize}
  These objects are required to satisfy the following conditions.  
  \begin{enumerate}
  \item
    For every $v\in V-\{v_0\}$ and
    $\tau\in C$, we have
    \[
         \tau^*\omega_{\tau(v)} = \chi(\tau)\cdot\omega_v.
    \]
  \item
    Let $v\in V$. We have $\delta\neq 0$ if and only if $v\neq
    v_0$. Moreover, if $\delta_v=1$ (resp.\ $0<\delta_v<1$) the
    differential $\omega_v$ is logarithmic (resp.\ exact). (If this
    holds, then we call the vertex $v$ {\sl multiplicative} (resp.\
    {\sl additive}.)
  \item
    For every edge $e\in E-\{e_0\}$, we have the equality
    \[
       \ord_{z_e}\omega_{s(e)} = -\ord_{z_e}\omega_{t(e)}-2.
    \]
  \item
    For every edge $e\in E$, we have
    \[
      \delta_{t(e)} = \delta_{s(t)}  + 
                      (p-1)\cdot\epsilon_e\cdot h_e,
    \]
    where 
    \[
        h_e  :=  \ord_{z_e}\omega_{t(e)}+1.
    \]
  \item
    For $b\in B$, let $Z_v$ be the component containing the point
    $z_b$. Then the differential $\omega_v$ has a simple pole in
    $z_b$. 
  \end{enumerate}
  The integer $h:=h_{e_0}$ is called the {\sl conductor} of the
  Hurwitz tree. The rational number $\delta:=\delta_{v_0}$ is called
  the {\sl different}. The {\sl tame inertia character} is the
  character $\lambda=\lambda_{e_0}:C\to k^\times$ which describes the
  action of $C$ on the tangent space of $Z$ at the point $z_0$. The
  family $(a_b)_{b\in B}$, with
  \[
     a_b  :=  \res_{z_b}\omega_v,
  \]    
  is called the {\sl type}.
\end{defn}

\begin{lem} \label{Htreelem1}
  Let $(Z,\omega_v,\delta_v,\epsilon_e)$ be a Hurwitz tree.
  \begin{enumerate}
  \item
    Fix an edge $e\in E$ and let $Z_e\subset Z$ be the union of all
    components $Z_v$ corresponding to vertices $v$ which are separated
    from the root $v_0$ by the edge $e$. Then
    \[
        h_e = |\{\,b\in B  \mid  z_b\in Z_e\,\}| \,-\, 1  >  0.
    \]
    In particular, $h=|B|-1 > 0$.
  \item
    The following conditions on a vertex $v$ are equivalent:
    \begin{enumerate}
    \item 
       $v$ is multiplicative,
    \item
       $v$ is a leaf of the tree $T$,
    \item
       the component $Z_v$ contains one of the marked points $z_b$.
    \end{enumerate}
  \item
    Recall that $\lambda_e:C_e\to k^\times$ is the character
    which describes the action of the stabilizer of the point $z_e$ on
    the cotangent space of $Z_{s(e)}$ at $z_e$. We have that 
    \[
          \lambda_e^{-h_e} = \chi|_{C_e}.
    \]
  \item
    For $b\in B$ and $\tau\in C$ we have
    $a_{\tau(b)}={\chi(\tau)}a_b$. 
  \item
    The character $\chi$ is either trivial or injective. If $\chi$ is
    injective then $m|h+1$. 
  \end{enumerate}
\end{lem}

\proof Statements (i) and (ii) are proved in \cite{GM99} and
\cite{Henrio02}. Statement (iii) is proved exactly as in the proof of
Proposition \ref{boundprop2}. Statement (iv) is an immediate
consequence of Condition (i) of Definition \ref{Htreedef}. It remains
to prove (v).

Write $m=m_1m_2$, where $m_1$ is the order of $\chi$ and $m_2$ is the
order of the kernel of $\chi$. Assume that $m_2\neq 1$. We want to
show that this implies $m_1=1$. Let $Z_v$ be the unique component of
$Z$ which contains the point $z_0$. By definition, the action of $C$
fixes the point $z_0$. Since $Z_v\cong\PP^1_k$, there exists a unique
point $z_1\not= z_0$ on $Z_v$ with nontrivial stabilizer, and in fact
$z_1$ is stabilized by the whole group $C$. It follows from the same
computation that we used to prove (iii) that
$m_2|\ord_{z_1}\omega_v+1$. In particular, $\ord_{z_1}\omega_v\not=0$.
Therefore, $z_1$ is either equal to a marked point $z_b$ or to a
singular point $z_e$. In the first case, it follows immediately from
 (iv) that $m_1=1$, as desired. In the second case, we
consider the subtree $Z_e\subset Z$ as in (i). It is easy to
check that $Z_e$ inherits from $Z$ the structure of a Hurwitz tree
of type $(C,\chi)$, with distinguished point $z_1$. Now the claim
$m_1=1$ follows by induction on the number of components of the
Hurwitz tree. This proves the first part of (v). Now assume that
$m_2=1$. Then (iv) implies that the action of $C$ on the set $B$ is
free. Therefore, $m|h+1$ follows from (i). This proves the second part
of (v), and completes the proof of the lemma.
\Endproof

\subsection{The Hurwitz tree associated to a $G$-action on the disc}  
  \label{Htree2}

Fix a cyclic group $C$ of order $m$, with $(p,m)=1$, and a character
$\chi:C\to\FF_p^\times$. We let $P$ denote a cyclic group of order
$p$, with generator $\sigma$. We set $G:=P\rtimes_\chi C$.
 
Let $R$ be a complete discrete valuation ring with residue field $k$.
Let $K$ denote the fraction field of $R$. We assume that $K$ has
characteristic $0$ and contains the $p$th roots of unity. Let
$Y_K:=\{\,y\,\mid\, |y|_K<1\,\}$ be the rigid open unit disc over $K$.
Suppose we are given a faithful action of $G$ on $Y_K$,
\[
   \phi_K:G  \inj  \Aut_K(Y_K).
\]
Obviously, $\phi_K$ induces a $G$-action on the boundary of the disc.
We let $h$ be the conductor, $\delta$ the different and $\lambda$ the
tame inertia character associated to this action (Definition
\ref{bounddef}). We assume that $h>0$.

Following \cite{Henrio02}, we will now associate to $\phi_K$ a Hurwitz
tree of type $(C,\chi)$.

Let $y_{b,K}\in Y_K$ be the fixed points of the automorphism
$\phi_K(\sigma)$, indexed by the finite set $B$. We assume that the
points $y_{b,K}$ are all $K$-rational. It is proved in \cite{GM99}
that $|B|=h+1$. Hence $|B|\geq 2$.  

Let $Y_R$ be the minimal model of $Y_K$ which separates the points
$y_{b,K}$. More precisely, $Y_R$ is an admissible blow-up of the formal
unit disc $\Spf R[[y]]$ such that 
\begin{itemize}
\item
  the exceptional divisor $Y$ of the blow-up $Y_R\to\Spf R[[y]]$ is a
  semistable curve over $k$,
\item
  the fixed points $y_{b,K}$ specialize to pairwise distinct smooth points
  $y_b$ on $Y$, and
\item
  if $y_0$ denotes the unique point on $Y$ which lies in the closure of
  $Y_R\otimes k-Y$, then $(Y,(y_b),y_0)$ is stably marked.
\end{itemize}
We call $(Y,(y_b),y_0)$ the {\sl special fiber} of the model $Y_R$. Note
that it is a decorated tree, in the sense of Definition
\ref{dectreedef}. Note also that there is a natural action of the
group $G$ on $(Y,(y_b),y_0)$. The element $\sigma$ fixes the points
$(y_b)$ and $y_0$. Therefore, $\sigma$ acts trivially on $Y$,
and we obtain a natural action of $C=G/P$ on $Y$. 

Let $Z_K:=Y_K/P$ be the quotient of the disc under the cyclic group of
order $p$; this is again a rigid open disc. Let $z_{b,K}\in Z_K$
denote the image of the fixed point $y_{b,K}$. Similarly to what we
did above, we define $Z_R$ as the minimal model of $Z_K$ which
separates the points $z_{b,K}$, and we let $(Z,(z_b),z_0)$ denote the
special fiber of $Z_R$. The canonical map $Y_K\to Z_K$ extends
uniquely to a map $Y_R\to Z_R$. One shows that the induced map $Y\to
Z$ is a $C$-equivariant homeomorphism, purely inseparable of degree
$p$.

Let $(V,E)$ be the combinatorial tree underlying 
$(Z,(z_b),z_0)$. (We will use freely the notation introduced in \S
\ref{Htree1}.) For $v\in V$, let $U_v\subset Z_v$ be the complement of
the singular and marked points and let $U_{v,K}\subset Z_K$ be the
affinoid subdomain with reduction $U_v$. Let $V_v\subset Y$ (resp.\
$V_{v,K}\subset Y_K$) denote the inverse image of $U_v$ (resp.\ of
$U_{v,K}$). By construction, $V_{v,K}\to U_{v,K}$ is a torsor under
the constant $K$-group scheme $P$. One shows that this torsor extends
to a torsor $V_{v,R}\to U_{v,R}$ under a certain finite flat $R$-group
scheme $P_v$ of degree $p$. Therefore, the finite inseparable map
$V_v\to U_v$ is naturally endowed with the structure of a torsor under
the finite, flat and local $k$-group scheme $P_v\otimes k$ of degree
$p$. See \cite{Henrio02} for details. According to the
classification of such group schemes, we have to distinguish two cases
(compare to the additive and the multiplicative case in \S \ref{bound2}).
\begin{itemize}
\item[(a)]
  $P_v\otimes k\cong\mu_{p,k}$: the $\mu_p$-torsor $V_v\to U_v$ can be
  classified by a logarithmic differential form $\omega_v$ on
  $U_v$. We set $\delta_v:=1$; this is the
  exponential valuation of the different of the
  $R$-group scheme $P_v\cong\mu_{p,R}$.
\item[(b)] $P_v\otimes k\cong\alpha_p$: the $\alpha_p$-torsor $V_v\to
  U_v$ can be classified by an exact differential form $\omega_v$ on
  $U_v$. We define $0<\delta_v<1$ as the exponential
  valuation of the different of the $R$-group
  scheme $P_v$.
\end{itemize} 
Finally, for $e\in E$ we let $A_e\subset Y_K$ denote the
subset of all points which specialize to the singular point $y_e\in Y$
corresponding to $e$. This is an open annulus. We define $\epsilon_e$
as the thickness of $A_e$, i.e.\ the positive rational number such
that
\[
     A_e  \cong   \{\,y  \mid  |p|_K^{\epsilon_e}<|y|_K<1\,\}.
\]

\begin{prop} \label{Htreeprop1}
  The datum $(Z,\omega_v,\delta_v,\epsilon_e)$ defines a Hurwitz tree
  of type $(C,\chi)$. Moreover, $h$ is the conductor, $\delta$ the
  different and $\lambda$ the tame inertia character of this Hurwitz tree.
\end{prop}

\proof Except for Condition (i) of Definition \ref{Htreedef} and the
statement on the tame inertia character which involve the $C$-action,
this is already proved in \cite{Henrio02}. To check Condition (i),
one easily verifies the following two facts: (a) the construction of
the datum $(Z,\omega_v,\delta_v,\epsilon_e)$ depends functorially on
the pair $(Y_K,\sigma)$, and (b) for $a\in\FF_p^\times$, the pair
$(Y_K,\sigma^a)$ gives rise to the datum
$(Z,a\cdot\omega_v,\delta_v,\epsilon_e)$.  The statement on the tame 
inertia character follows from the fact that the action of $C$ on the
semistable model $Z_R$ is admissible. \Endproof

\begin{cor} \label{Htreecor1}
  Let $\phi:G\inj\Aut_k(k[[y]])$ be a local $G$-action. If $\phi$
  lifts to characteristic zero, then the character $\chi$ is either
  trivial or injective. Moreover, if $\chi$ is injective then $m|h+1$.
\end{cor}

\proof
A lift of $\phi$ induces a $G$-action on the open rigid disc (with
discriminant $0$). Hence the corollary follows from Proposition
\ref{Htreeprop1} and Lemma \ref{Htreelem1} (v).
\Endproof

\subsection{Construction of $G$-actions on the disc} \label{Htree3}

In \S \ref{Htree2} we have associated to a $G$-action
$\phi_K:G\inj\Aut_K(Y_K)$ on the open unit disc a Hurwitz tree of type
$(C,\chi)$. The main result of this section is that this construction
can be reversed.

\begin{thm} \label{Htreethm}
  Every Hurwitz tree $(Z,\omega_v,\delta_v,\epsilon_e)$ of type
  $(C,\chi)$ is associated to a faithful $G$-action
  $\phi_K:G\inj\Aut_K(Y_K)$ on the open rigid unit disc over $K$, for
  some finite extension $K$ of ${\rm Frac}(W(k))$. 
\end{thm}

The proof of Theorem \ref{Htreethm} will be given in \S \ref{Htree6}.
The relevance for the local lifting problem is summarized in the next
corollary.

\begin{cor} \label{Htreecor}
  Let $\phi:G\inj\Aut_k(k[[y]])$ be a local $G$-action with conductor
  $h$. Then $\phi$ lifts to characteristic $0$ if and only if there
  exists a Hurwitz tree of type $(C,\chi)$, conductor $h$ and
  discriminant $\delta=0$.
\end{cor}

\proof The `only if' part of the corollary follows already from
Proposition \ref{Htreeprop1}. Let $Z$ be a Hurwitz tree of type
$(C,\chi)$, conductor $h$ and discriminant $\delta=0$. By Theorem
\ref{Htreethm}, there exists a $G$-action on the open rigid disc which
gives rise to the Hurwitz tree $Z$. Let $\phi_R:G\inj\Aut_R(R[[y]])$
be the induced action on the formal model $\Spf R[[y]]$ of $Y_K$. The
reduction of $\phi_R$ to $k$ corresponds to a local $G$-action
$\phi':G\inj\Aut_k(k[[y]])$ with conductor $h$. Let $\lambda:C\to
k^\times$ be the tame inertia character associated to $\phi_K$. Up to
changing $\phi$ by an automorphism of $G$, we may assume that
$\phi(\tau)^*y\equiv \lambda(\tau)\cdot y\pmod{y^2}$.  Remark
\ref{boundrem1} implies that $\phi$ and $\phi'$ are conjugate in
$\Aut_k(k[[y]])$. Hence we may regard $\phi_R$ as a lift of $\phi$.
\Endproof

\subsection{Realization of the vertices} \label{Htree4}

Let $(C,\chi)$ be as before, and set $G:=P\rtimes_\chi C$. Set
$Z:=\PP^1_k$. We let $C$ act on $Z$ in such a way that
$\tau^*z=\lambda(\tau)\cdot z$ for all $\tau$, where $z$ is the
standard parameter on $\PP^1$ and $\lambda:C\to k^\times$ is a
character of order $m=|C|$. Suppose we are given a differential form
$\omega$ on $Z$ and a rational number $0<\delta\leq 1$ such that the
following holds.

\begin{itemize}
\item The differential $\omega$ has poles in $r\geq 2$ points
  $z_1,\ldots,z_r\not=\infty$, a zero at $\infty$ and no other poles
  or zeroes. 
\item
  For all $\tau\in C$ we have $\tau^*\omega=\chi(\tau)\cdot\omega$.
\item
  If $\delta=1$ (resp.\ $\delta<1$) then $\omega$ is logarithmic
  (resp.\ exact).  
\end{itemize}
Note that the datum $(Z,\omega,\delta)$, together with the $C$-action,
essentially corresponds to a vertex $v\not=v_0$ of a Hurwitz tree of
type $(C,\chi)$. The points $z_i$ (resp.\ the point $\infty$)
correspond to the singular points $z_e$ such that $s(e)=v$ (resp.\ to
the unique singular point $z_e$ with $t(e)=v$).

Suppose first that $\delta=1$ (the multiplicative case). Then we set
$R:=W(k)[\zeta_p]$ and $K:={\rm Frac}(R)$. We fix an isomorphism
$P\cong\mu_p(K)$; this choice allows us to view the abstract group $G$
as a group scheme over $K$.  Then $G_R:=\mu_{p,R}\rtimes_\chi C$ is
the unique finite flat group scheme over $R$ with generic fiber $G$
and with connected subgroup scheme $P_R:=\mu_{p,R}$. We write
$G_k:=G_R\otimes_R k$.  We set $U:=Z-\{\infty\}=\AA^1_k$.  Then the
logarithmic differential $\omega$ gives rise to a finite flat and
radicial morphism of smooth curves $V\to U$, together with an action
of $P_k=\mu_{p,k}$ on $V$ such that $U=V/P_k$. The formula
$\tau^*\omega=\chi(\tau)\cdot\omega$ shows that this action extends
uniquely to an action of $G_k$ on $V$ which induces the canonical
action of $C=G_k/P_k$ on $U$. See \cite[\S 4.1]{cotang} for details.

Now suppose that $0<\delta<1$ (the additive case). Write
$\delta=1-n/d$, for some $n,d\in\NN$.  Let $K$ be an extension of
${\rm Frac}(W(k)[\zeta_p])$ of degree $d$ and let $\pi$ denote a
uniformizer of $K$. Let $\Hc_n$ be the finite flat group scheme over
$R$ defined in \cite[\S 1.1]{Henrio02}. We have $\Hc_n\otimes
K\cong\mu_{p,K}$ and $\Hc_n\otimes k\cong\alpha_p$.  Similar to what
we did in the multiplicative case, we identify $\Hc_n(K)$ with $P$ and
let $G_R:=\Hc_n\rtimes_\chi C$ denote the unique finite flat group
scheme over $R$ with generic fiber $G$ and with connected subgroup
scheme $P_R:=\Hc_n$. We write $G_k:=G_R\otimes_R k$ and set
$U:=Z-\{z_1,\ldots,z_r,\infty\}$. Then the exact differential $\omega$
gives rise to a $P_k=\alpha_p$-torsor $V\to U$. Moreover, the
$P_k$-action extends uniquely to an action of $G_k$ which induces the
given action of $C=G_k/P_k$ on $U$.

\begin{lem}
  The curve $V$ lifts to a smooth formal scheme $V_R$ over $R$,
  together with an action of $G_R$ which lifts the action of $G_k$ on
  $V$. 
\end{lem}

\proof We give an abstract proof of this fact which treats the
multiplicative case and the additive case simultaneously. It is not
hard to give a more down-to-earth proof, using the explicit equations for
$\mu_p$- and $\Hc_n$-torsors. 

Let $\Def(V,G_R)$ be the functor which classifies $G_R$-equivariant
deformations of $V$ over local artinian $R$-algebras with residue
field $k$. To prove the lemma, it suffices to show that $\Def(V,G_R)$
is unobstructed. By \cite{BM00} and \cite{cotang}, the
obstructions for $\Def(V,G_R)$ are represented by elements of the
equivariant cohomology group $H^2(G_k,V,\T_{V/k})$ (here $\T_{V/k}$
denotes the tangent sheaf on $V$; it affords a natural action of
$G_k$). The spectral sequence
\[
     E_2^{i,j}=H^i(C,H^j(P_k,V,\T_{V/k}))  \Rightarrow  
          H^{i+j}(G_k,V,\T_{V/k}),
\]
together with the fact that the order of $C$ is prime to $p$, shows
that 
\[
      H^2(G_k,V,\T_{V/k}) = H^2(P_k,V,\T_{V/k})^C. 
\]
Moreover, if an element $c$ of $H^2(G_k,V,\T_{V/k})$ represents an
obstruction for the functor $\Def(V,G_R)$, then its image in
$H^2(P_k,V,\T_{V/k})$ represents the induced obstruction for the
functor $\Def(V,P_R)$. Hence it suffices to show that $\Def(V,P_R)$ is
unobstructed. 

Suppose we are in the multiplicative case, and let $A$ be a local
artinian $R$-algebra with residue field $k$. Let $V_A$ be a
$P_R$-equivariant lift of $V$. Since $V$ is affine, $V_A$ is affine as
well. Therefore, the quotient map $V_A\to U_A:=V_A/P_R$ is given by a
global Kummer equation $y^p=f$. If $A\to A'$ is a small extension,
we can construct a $G_R$-equivariant lift of $V_A$ to $A'$ by first
lifting $U_A$ and then lifting the Kummer equation. This shows that
$\Def(V,P_R)$ is unobstructed. The proof in the additive case is
essentially the same, using the explicit equations for $\Hc_n$-torsors
given in \cite[\S 1.1]{Henrio02}. This completes the proof of the
lemma.
\Endproof

Let $V_R$ be a lift of $V$ over $R$, together with an action of $G_R$
lifting the action of $G_k$ on $V$. Set $V_K:=V_R\otimes_R K$. This is
an affinoid with reduction $V$ and carries an action of the abstract
group $G=G_K$. 

In the multiplicative case we have $U=\AA^1_k$ and hence
$V\cong\AA^1_k$. Therefore, $V_K\cong\{\,y\,\mid\, |y|\leq 1\,\}$ is
isomorphic to a closed unit disc.

In the additive case we have $U=\AA^1_k-\{z_1,\ldots,z_r\}$ and hence
$V\cong\AA^1_k-\{y_1,\ldots,y_r\}$. Therefore,
\[
    V_K  \cong  \{\,y \mid  |y|\leq 1,\, |y-y_{i,K}|\geq 1 \}
\]
is isomorphic to the complement of $r$ open discs inside one closed
disc. Let $\Spec S_i$ (resp.\ $\Spec S_\infty$) be the boundary of the
missing open disc containing the point $y_{i,K}$ (resp.\ $\infty$). In
other words, $S_i=R[[y-y_{i,K}]]\{(y-y_{i,K})^{-1}\}$ and
$S_\infty=R[[y^{-1}]]\{y\}$. Let $G_i:=P\rtimes C_i$ be the stabilizer
of the point $y_i$. By construction, we have the following result. 

\begin{prop} \label{vertexprop}
  The action of $G$ on $V_K$ induces an action of $G_i$ on $\Spec S_i$
  (resp.\ of $G$ on $\Spec S_\infty$). This action has conductor
  $h_i:=\ord_{z_i}\omega+1$ (resp.\ conductor
  $h:=-\ord_\infty\omega-1$) and discriminant $\delta$.
\end{prop}

\subsection{Realization of the edges} \label{Htree5}

We let $G=P\rtimes_\chi C$ be as before. Suppose we are given the
following data. 
\begin{itemize}
\item
  An integer $h\geq 1$ with $\ord(\chi)|h+1$ and $(h,p)=1$. 
\item
  An injective character $\lambda:C\inj k^\times$ with
  $\lambda^{-h}=\chi$. 
\item
  Rational numbers $\delta_0,\delta_1$ with
  $0\leq\delta_0<\delta_1\leq 1$. 
\end{itemize}
These data essentially correspond to an edge of a Hurwitz
tree of type $(C,\chi)$. We set 
\[
   \epsilon  :=  \frac{\delta_1-\delta_0}{(p-1)h}, \qquad
   \kappa_i  :=  \frac{1-\delta_i}{(p-1)h},  \;i=0,1.
\]
Let $K$ be a sufficiently large finite and tamely ramified extension
of the fraction field of $W(k)$. Let $R$ denote the ring of
integers of $K$. Choose a $p$th root of unity $\zeta_p\in R$.
Let $\rho_i\in R$ be an element with $|\rho_i|=|p|^{\kappa_i}$, and
set $\pi:=\rho_0/\rho_1$. Then
\[
    B  :=  \{\, z_0 \,\mid\, |\pi^{p}| < |z_0| < 1 \,\}
\]
is an open annulus of thickness $p\epsilon$, with formal model $\Spf
R[[z_0,z_1\mid z_0z_1=\pi^p]]$. We let the group $C$ act on
$B$ such that $\tau^*z_0=\lambda(\tau)\cdot z_0$. 

Choose a generator $\tau$ of $C$ and a generator $\sigma$ of $P$. Set
$\zeta_m:=\lambda(\tau)$.  Write $m=m_1m_2$, where $m_1|(p-1)$ is the
order of the character $\chi$. Choose an integer $\nu$ with
$\nu\equiv\zeta^{-h}\pmod{p}$. Then $1-\nu^{m_1}=np$, for some integer
$n$. We define
\begin{equation} \label{edgeeq1}
   f  :=  \prod_{i=0}^{m_1-1} (1+\zeta_m^{ih}\rho_1^{ph}z_1^h)^{\nu^i}.
\end{equation}
Clearly, $f$ is a bounded invertible function on
$B$. A straightforward computation shows that
\begin{equation} \label{edgeeq2}
    \tau^*f = f^\nu\cdot(1+\rho_1^{ph}z_1^h)^{np}.
\end{equation}
Let $A\to B$ be the $\mu_p$-torsor defined by the Kummer equation
$w^p=f$. One checks that the $C$-action on $B$ extends to a $G$-action
on $A$ such that
\begin{equation} \label{edgeeq3}
   \sigma^*w = \zeta_p\cdot w, \qquad
   \tau^*w = w^\nu\,(1+\rho_1^{ph}z_1^h)^n.
\end{equation}

\begin{prop}  \label{edgeprop}
  \begin{enumerate}
  \item
    The rigid space $A$ is an open annulus of thickness $\epsilon$,
    i.e.\ 
    \[
         A \cong \{\, y \,\mid\, |\pi|<|y|<1 \,\}.
    \]
  \item
    Let $\Spec S_0=\Spec R[[y]]\{y^{-1}\}$ be the `outer boundary' of $A$
    (which is mapped to $\Spec R[[z_0]]\{z_0^{-1}\}$). The induced action
    of $G$ on $\Spec S_0$ has conductor $h$, discriminant $\delta_0$ and
    tame inertia character $\lambda$.
  \item Let $\Spec S_1=\Spec R[[\pi y^{-1}]]\{\pi y\}$ be the `inner
    boundary' of $A$ (which is mapped to $\Spec
    R[[z_1]]\{z_1^{-1}\}$). The induced action of $G$ on $\Spec S_1$ has
    conductor $-h$, discriminant $\delta_1$ and tame inertia character
    $\lambda^{-1}$.
  \end{enumerate}
\end{prop}

\proof 
Equation \eqref{edgeeq1} shows that
\[
  f = 1 + 
       (\sum_{i=0}^{m_1-1}(\zeta_m^h\nu)^i)\rho_0^{ph}\,z_0^{-h}+\ldots
    = 1 + m_1\rho_0^{ph}\,z_0^{-h}+\ldots
\]
(we have used that $\zeta_m^h\nu\equiv 1\pmod{p}$). The proposition
follows now from the explicit computations in \cite[\S 1]{Henrio02}. 
\Endproof

\subsection{The proof of Theorem \ref{Htreethm}} \label{Htree6}

Let $Z=(Z,\omega_v,\delta_v,\epsilon_e)$ be a Hurwitz tree of type
$(C,\chi)$, conductor $h$, discriminant $\delta$ and tame inertia
character $\lambda$. We call an action of the group $G=P\rtimes_\chi
C$ on the open disc a {\sl realization} of $Z$ if $Z$ is associated to
this action by the construction of \S \ref{Htree2}. Theorem
\ref{Htreethm} claims that we can realize $Z$. We will prove this
claim by induction on the depth of the tree $Z$. Therefore, we may
assume that we can realize all Hurwitz trees (of some type
$(C',\chi')$) with lower depth than $Z$.

Let $e_0$ be the edge whose source is the root $v_0$, and set
$v_1=t(e_0)$. Then $Z_0:=Z_{v_1}$ is the unique component of $Z$ which
contains the distinguished point $z_0\in Z$. Let
$\delta_1:=\delta_{v_1}$. Let $e_1,\ldots,e_r$ be the edges with
source $v$ and $z_i:=z_{e_i}$ the corresponding singular points. Set
$h_i:=h_{e_i}$. Let $C_i\subset C$ denote the stabilizer of the point
$z_i$, and set $G_i:=P\cdot C_i\subset G$. Note that either $C_i=C$ or
$C_i=1$. Let $Z_i\subset Z$ be the subtree of $Z$ which contains the
point $z_i$ but not the component $Z_0$. It is clear that $Z_i$
inherits from $Z$ the structure of a Hurwitz tree of type
$(C_i,\chi|_{C_i})$, with conductor $h_i$ and different $\delta_1$. By
our induction hypothesis, there exists an open unit disc $Y_{i,K}$
over some finite extension $K$ of ${\rm Frac}(W(k))$, together with an
action of $G_i$, whose associated Hurwitz tree is $Z_i$. The extension
$K$ may be taken to be independent of $i$. We may
also assume that the discs $Y_{i,K}$ corresponding to one $G$-orbit of
the trees $Z_i$, together with the group actions, are all isomorphic.
Then we can define an action of $G$ on the disjoint union of the discs
$Y_{i,K}$ which induces the $G_i$-action on $Y_{i,K}$. Let $\Spec S_i$
be the boundary of $Y_{i,K}$.

By Lemma \ref{Htreelem1}, the differential form $\omega:=\omega_{v_1}$
on $Z_0$ has a zero of order $h-1$ at $z_0$ and poles of order $h_i+1$
at $z_i$. Let $V_K$ be the affinoid with $G$-action constructed in \S
\ref{Htree4}, starting from the datum $(Z_0,\omega,\delta_1)$. By
Proposition \ref{boundprop3} and Proposition \ref{vertexprop}, we can
identify the boundary of the missing open disc corresponding to the
point $z_i$ with $\Spec S_i$, in a way which is compatible with the
action of $G_i$. Clearly, we can choose this identification such that
the $G$-action on the disjoint union of the boundaries $\Spec S_i$
induced from the $G$-action on the disjoint union of the discs
$Y_{i,K}$ agrees with the $G$-action on $V_K$. We can now use
\cite[Lemme 3.7]{Henrio02} to patch together the affinoid $V_K$ and
the discs $Y_{i,K}$, in a $G$-equivariant way.  The result is a closed
disc $W_K\cong\{\,w\,\mid\,|w|\leq 1\,\}$, together with an action of
$G$. By Proposition \ref{vertexprop} and the construction, the induced
action on the boundary of $\Spec R[[w^{-1}]]\{w\}$ of $W_K$ has
conductor $h^{-1}$, discriminant $\delta_1$ and tame inertia character
$\lambda^{-1}$.

Let $A_K$ be the open annulus with $G$-action constructed in \S
\ref{Htree5}, starting from the datum
$(h,\lambda,\delta_0:=\delta,\delta_1)$. By Proposition
\ref{boundprop3} and Proposition \ref{edgeprop}, we can identify the
boundary of $W_K$ with the `inner boundary' of $A_K$, in a
$G$-equivariant way. By \cite[Lemme 3.8]{Henrio02} we can patch
together $W_K$ and $A_K$, together with the $G$-action, along these
boundaries. The result is an open disc $Y_K$ with $G$-action, of
conductor $h$, discriminant $\delta$ and inertia character $\lambda$.
By construction, $Y_K$ is a realization of the Hurwitz tree $Z$. This
completes the proof of Theorem \ref{Htreethm}.
\Endproof

\begin{rem}
  \begin{enumerate}
  \item
    It would be more satisfactory (and this paper would be shorter) if one
    could prove Theorem \ref{Htreethm} using the results of
    \cite{Henrio02} (i.e.\ the case $G=\ZZ/p$) as a `black
    box'. However, since the automorphism of the
    disc of order $p$ constructed in \cite{Henrio02} does not depend
    in a functorial way on the Hurwitz tree, this looks difficult. 
  \item  
    Let $\sigma:Y_K\iso Y_K$ be an automorphism of order $p$ of the
    disc with conductor $h<p$. Then \cite[Theorem III.3.1]{GM99} shows
    that the Hurwitz tree $Z$ associated to $\sigma$ is irreducible,
    i.e.\ consists of a single component $Z\cong\PP^1$. Furthermore,
    \cite[Theorem V.6.3.1]{GM99} shows that the conjugacy class of
    $\sigma$ in $\Aut_K(Y_K)$ is uniquely determined by the Hurwitz
    tree $Z$. If this would be true for $h>p$ and for arbitrary Hurwitz
    trees then the proof of Theorem \ref{Htreethm} could be simplified
    considerably, as suggested in (i).
  \end{enumerate}
\end{rem}


\section{Local lifts for large conductor}\label{largehsec}

\subsection{Existence of differential forms}\label{existencesec}
In constructing a Hurwitz tree as in \S \ref{Htree1}, the main
difficulty is to find the differential forms $\omega_v$. In this
section, we illustrate the problems in the easiest case, namely the
case that $G$ is cyclic. We will use these results afterwards to prove
the local lifting result in case $h>p$.

\begin{lem}\label{omegalem}
  Let $h>0$ be an integer which is prime to $p$. There exists a
  logarithmic differential on $\PP^1_k$ which has $h+1$ simple poles
  and a single zero of order $h-1$.
\end{lem}

\proof Choose a primitive $h$th root of unity $\zeta_h\in k$ and
define $z_i=\zeta_h^i\in \PP^1_k$ for $i=1, \dots, h$. Put
$z_0=0$. Let $u=z^{-h}\prod_{i=1}^h(z-z_i)$. Then
\[
  \omega:=\frac{{\rm d}u}{u}=\left(\frac{-h}{z}+\sum_{i=1}^h
  \frac{1}{z-z_i}\right){\rm d}z=
\frac{h}{z^{h+1}-z}{\rm d}z.
\]
Note that $\omega$ has a single zero at $z=\infty$.   \Endproof

A  result in the same direction is proved in \cite[\S
  3.5]{Henrio02}.

We now discuss several problems one encounters while explicitly
constructing Hurwitz trees.  If $h+1\leq p$ then the Hurwitz tree $Z$
is necessarily irreducible (\cite{GM99}). As soon as $h+1>p$ one
has greater freedom in choosing the tree $Z$. This makes the problem
easier.

Suppose that the Hurwitz tree $Z$ is irreducible, and let
$\mathbf{a}=(a_b)$ be the type of $Z$ (see Definition \ref{Htreedef}).
In general it is not easy to decide which possibilities for
$\mathbf{a}$ occur. Lemma \ref{omegalem} may be rephrased as follows.
There exists a logarithmic differential $\omega$ on $\PP^1_k$ with
residues ${\mathbf a}=(1,1,\ldots ,1,p-h)$ and a single zero.  The
following lemma shows that not all sets of integers $\mathbf a$ whose
sum is zero in $\FF_p$ occur as set of residues of a logarithmic
differential.

\begin{lem} \label{omegalem2}
Let $p=5$ and ${\boldsymbol a}=(1,1,4,4)$. There exists no logarithmic
differential on $\PP^1_k$ with a single zero whose set of residues is
$\boldsymbol a$.
\end{lem}

\proof Suppose that $\omega$ is a logarithmic differential on
$\PP^1_k$ with a single zero at $z=\infty$ and set of residues
${\boldsymbol a}=(1,1,4,4)$. We may suppose that $\omega$ has poles in
the four pairwise distinct points $z=0,1,\lambda, \mu$. It follows
that $\omega$ is the logarithmic differential of $
f=uz(z-1)(z-\lambda)^4(z-\mu)^4,$ for some $u\in k^\ast$. We find that
\[
   \omega=\frac{1}{z}+\frac{1}{z-1}+\frac{-1}{z-\lambda}+\frac{-1}{z-\mu}
         = \frac{(1-\lambda-\mu)z^2+2\lambda\mu
         z-\lambda\mu}{z(z-1)(z-\lambda)(z-\mu)}.
\]

Since $\omega$ has a single zero at $z=\infty$, we should choose
$\lambda, \mu$ such that $\lambda+\mu=1$ and $\lambda\mu=0$. But this
contradicts the hypothesis that $0,1,\lambda, \mu$ are pairwise
distinct.
\Endproof

The reason the differentials of Lemma \ref{omegalem} are so easy to
write down is that there is an extra automorphism (of order $h$) which
fixes the differential. If $G=D_p$ this trick no longer works, and it
will be much harder to find a set of residues ${\boldsymbol a}$ for
which we can find a differential.

\begin{rem} \label{Leorem}
  Leonardo Zapponi has proved a general criterion for the existence of
  logarithmic differentials with a single zero which includes Lemma
  \ref{omegalem} and Lemma \ref{omegalem2} as special cases.
\end{rem}


\subsection{}

Let $p$ be an odd prime and $G$ the dihedral group of order $2p$. We
write $C$ for the quotient of $G$ of order two and let $\tau$ be its
generator. Let $\chi:C\iso\{\pm 1\}$ be the character of order $2$. 

\begin{thm}\label{largehthm}
Let $\phi:G\inj\Aut_k(k[[y]])$ be a local $G$-action. Let $h$ be the
conductor and suppose that $h>p$. Then $\phi$ lifts to characteristic
zero.
\end{thm}

\proof By Corollary \ref{Htreecor}, it suffices to construct a Hurwitz
tree $Z=Z_h$ of type $(C,\chi)$, conductor $h$ and discriminant $0$.
The Hasse--Arf theorem implies that the conductor $h$ is odd, and we
may write $h+1=2(\alpha+p\beta)$, with $(p+3)/2\leq \alpha\leq
(3p-1)/2$. Since $h$ is prime to $p$, we exclude $\alpha\equiv
(p+1)/2\bmod{p}$.

\bigskip\noindent {\bf Case 1}: $(p+3)/2\leq \alpha\leq p$.  We
construct a tree $Z$ as follows. The vertex set consists of $v_0, v_1,
w_j, u_j$, for $j=0,\ldots, \beta$. There is one
edge connecting $v_0$ with $v_1$, and for $j=0,\ldots,\beta$ there is
an edge $e_j$ (resp.\ $f_j$) connecting $v_1$ with $w_j$ (resp.\
$u_j$).

For each $v\in V$, we choose a curve $Z_v$ of genus zero together with
a coordinate which we denote by $z$. We suppose that the edge $e$ with
source $v_1$ and target $v_0$ corresponds to $z_e=\infty\in
Z_{v_1}$. We let $\tau$ act on $Z_{v_1}$ by $\tau^* z=-z$. Choose
points $z_{e_j}\in Z_{v_1}-\{0, \infty\}$ and define $z_{f_j}=\tau
( z_{e_j})=-z_{e_j}$. Choose the points so that all $2(\beta+1)$
points are pairwise distinct. It is no restriction to suppose that
$z_{e_0}=1$.

 We construct an exact differential $\omega_{v_1}$ on
$Z_{v_1}$ with a pole of order $\alpha$ in $z=\pm z_{e_0}=\pm 1$ and a
pole of order $p$ in $z=\pm z_{e_j}$ for $j\neq 0$ by
\[
  \omega=\frac{{\rm d} z}{(z^2-1)^\alpha\prod_{j=1}^\beta
  (z^2-z_{e_j}^2)^p}.
\]
It is easy to see that this differential is exact. Here we use the
assumption that $(p+3)/2\leq \alpha\leq p$.

On the component $Z_{w_j}$ we define a logarithmic differential with a
single zero in some point $z_{e_j}$ and $\alpha$ (resp.\ $p$) simple
poles if $j=0$ (resp.\ $j\neq 0$). The existence of such differentials
follows from Lemma \ref{omegalem}. We define the pair
$(Z_{u_j},\omega_{u_j})$ as a copy of the pair
$(Z_{w_j},\omega_{w_j})$. It is clear that the action of $C$ on
$Z_{v_1}$ extends to an admissible action on the whole tree $Z$ which
verifies Condition (i) of Definition \ref{Htreedef}.

It remains to choose the discriminants $\delta_v$ and the thicknesses
$\epsilon_e$. We set $\delta_{v_0}:=0$,
$\delta_{w_j},\delta_{u_j}:=1$, and for $\delta_{v_1}$ we may choose
any rational number with $0<\delta_{v_1}<1$. Then Condition (iv) of
Definition \ref{Htreedef} imposes a unique value for all $\epsilon_e$.
This completes the construction of the Hurwitz tree $Z$ and the proof
of the theorem in Case 1.

\bigskip\noindent {\bf Case 2}: $p+1\leq \alpha \leq (3p-1)/2$.  Write
$ \alpha=\alpha_1+\alpha_2,$ with $ 1<\alpha_2\leq \alpha_1< p,$ note
that this is always possible.  We construct a tree $Z$ together with
an automorphism $\tau$ of order two as follows. The vertex set of $Z$
consists of vertices $v_0, v_1, w_{01}, u_{01}, w_{02}, u_{02}, w_1,
u_1,\ldots, u_{\beta}$. We choose components $Z_{v_0}, Z_{v_1},
Z_{w_{01}}, \ldots Z_{u_{\beta}}$, together with an action and the
automorphism $\tau$ of order two as in Case 1, i.e.\ $\tau$ sends
$w_j$ to $u_j$.  We denote again by $e_j$ (resp.\ $f_j$) the edge with
source $v_1$ and target $w_j$ (resp.\ $u_j$) and choose corresponding
points $z_{e_j}$ (resp.\ $z_{f_j}$) on $Z_{v_1}$, as follows. We take
$z_{e_{01}}=1$, $z_{f_{01}}=-1$, $z_{e_{02}}=\lambda$ and
$z_{f_{02}}=-\lambda$, where $\lambda$ will be specified later on in
the proof. For $j=1,\ldots, \beta$, we choose points $z_{e_j}\in
Z_{v_1}-\{0,\infty\}$ and define $z_{f_j}=-z_{e_j}$. We choose these
points in such a way that all points $\pm1, \pm \lambda, \pm z_{e_j}$
are pairwise distinct.

It follows from Lemma \ref{omegalem} and the fact that $\alpha_j< p$
that we may choose for  $j=1,2$ a logarithmic differential
$\omega_{w_{0j}}$ which has simple poles in $\alpha_j$ points on
$Z_{w_{0j}} -\{\infty\}$ and a single zero at $\infty$. We define
$\omega_{u_{0j}}=\tau^\ast \omega_{w_{0j}}$. Similarly, for 
$j=1, \ldots, \beta$, there exist logarithmic differentials
$\omega_{w_j}$ on $Z_{w_j}$ with $p$ poles and a single zero. 

On $Z_{v_1}$, we want to find a differential $\omega_{v_1}$ with  a pole
of order $\alpha_j$ in $\pm z_{e_{0j}}$ for $j=1,2$, and poles
of order $p$ in $\pm z_{e_j}$, for $j=1,\ldots, \beta$. To simplify the
notation, we write $z_j=z_{e_j}$. After
multiplying $\omega_{v_1}$ with an element of $k^\times$, the
differential is given by 
\[
   \omega_{v_1}=\frac{{\rm
   d}z}{(z^2-1)^{\alpha_1}(z^2-\lambda^2)^{\alpha_2}\prod_{j=1}^\beta
   (z^2-z_j^2)^p}=\frac{\eta\,{\rm
   d}z }{Q^p},
\]
where $\eta:=(z^2-1)^{p-\alpha_1}(z^2-\lambda^2)^{p-\alpha_2}$ and
$Q:=(z^2-1)(z^2-\lambda)\prod_j(z^2-z_j^2).$ We claim that
$\omega_{v_1}$ is exact, for suitable choice of $\lambda$.
Write $\eta=\sum_{i=0}^{2p-\alpha}\eta_i z^{2i}$.

\begin{lem}
If $\eta_{(p-1)/2}=0$, then the differential $\omega_{v_1}$ is exact.
\end{lem}

\proof
Note that the degree of $\eta$ in $x$ is $4p-2\alpha\leq 2p-2$. It
follows that the only $0\leq i\leq p-1$ for which $2i$ is congruent to
$-1\bmod{p}$ is $i=(p-1)/2$. Suppose that $\eta_{(p-1)/2}=0$ and
define
\[ 
   G=\frac{1}{Q^p}\sum_{i=0}^{2p-\alpha}\frac{\eta_i x^{2i+1}}{2i+1}.
\]
Then $\omega_{v_1}={\rm d}G$.
\Endproof

The lemma implies that we have to choose $\lambda$ in such a way that
$\eta_{(p-1)/2}=0$. One easily computes that 
\[
   \eta_{(p-1)/2}=\pm
   \sum_{i+j=(3p+1)/2-\alpha}\binom{p-\alpha_1}{i}\binom{p-\alpha_2}{j}
\lambda^j.
\]
The degree of $\eta_{(p-1)/2}$ in $\lambda$ is $\min(p-\alpha_2,
(3p+1)/2-\alpha)$. The polynomial has a zero of order $\max(0,
(p+1)/2-\alpha_2)$ in $\lambda=0$ and a zero of order $\max(0,
(p+1)/2-\alpha_2)$ at $\lambda=1$. (The first two statements are
obvious and the second one follows by symmetry.) 

Recall that $2\leq \alpha_2\leq \alpha_1\leq p-1$ and $p+1\leq
\alpha=\alpha_1+\alpha_2\leq (3p-1)/2$. Therefore $\alpha_1\geq (p+1)/2$. This
implies that the degree of $\eta_{(p-1)/2}$ is $(3p+1)/2-\alpha$ and
that $\eta_{(p-1)/2}$ has no zero at $\lambda=1$. It follows that the
number of zeroes of $\eta_{(p-1)/2}$ different from $0,1$ is 
\[
  \left\{\begin{array}{ll}
  (3p+1)/2-\alpha-(p+1)/2+\alpha_2=p-\alpha_1\geq 1&\mbox{ if }
  \alpha\leq (p+1)/2,\\ (3p+1)/2-\alpha\geq 1& \mbox{ if } \alpha\geq
  (p+1)/2.
   \end{array}\right.
\]
This implies that it is possible to choose $\lambda$ such that
$\eta_{(p-1)/2}=0$ and therefore such that $\omega_{v_1}$ is exact. The
proof of the proposition in this case follows now as in Case 1.
\Endproof


\section{Local lifts for small conductor}\label{smallhsec}
In this section we prove  the local lifting problem
for conductor $h<p$. Section \ref{trivsolsec} reformulates the problem
in terms of concrete equations, and describes the so called trivial
solutions of these equations. (These are the solutions that do not
corresponds to a solution of our problem.) In \S \ref{mainthmsec}
we use this to solve the lifting problem for small conductor. Section
\ref{largepsec} contains a weaker version of this result, namely we
suppose that $p$ is large with respect to $h$. 

\subsection{The trivial solutions}\label{trivsolsec}

Let $p$ be an odd prime and $G$ the dihedral group of order $2p$. Let
$C$ be its quotient of order two, $\tau$ the generator of $C$ and
$\chi:C\iso\{\pm 1\}$ the unique character of order $2$.  Let $0<h<p$
be an odd integer; write $\alpha=(h+1)/2$. To prove that all local
$G$-actions with conductor $h$ lift to characteristic zero, we need to
construct a Hurwitz tree $Z=Z_h$ of type $(C,\chi)$, conductor $h$ and
discriminant $0$, as in \S \ref{largehsec}. As explained in \S
\ref{existencesec}, this Hurwitz tree is irreducible, hence we may
suppose $Z=\PP^1_k$. We may also suppose that the distinguished point
$z_0=\infty$ and that $\tau$ acts on the standard parameter $z$ of $Z$
as $z\mapsto -z$. Our problem is thus reduced to finding a logarithmic
differential form $\omega$ on $\PP^1_k$ with $h+1=2\alpha$ poles and a
single zero at $\infty$ such that $\tau^*\omega=-\omega$.

Let ${\mathbf z}=(z_1, \ldots, z_\alpha)$ be points of $Z-\{0,
\infty\}$ such that $z_i^2\neq z_j^2$ if $i\neq j$. Suppose ${\mathbf
  a}=(a_1, \ldots, a_\alpha)\in (\FF_p^\times)^\alpha$ is a set of
residues which we consider to be fixed in this section. Put
\[
   f=f_{{\mathbf
z};{\mathbf a}}=\prod_{i=1}^\alpha\left(\frac{z-z_i}{z+z_i}\right)^{a_i},
\quad\mbox{ and }\quad\omega=\omega_{{\mathbf z};{\mathbf a}}=\frac{{\rm
d}f}{f}=\sum_{i=1}^\alpha \frac{2a_iz_i}{z^2-z_i^2}{\rm d}z.
\]

\begin{lem}\label{hsmalllem}
  Suppose that there exists a vector ${\mathbf z}$ as above such that
  $\omega_{{\mathbf z};{\mathbf a}}$ has a single zero at $\infty$.
  Then every local $G$-action with conductor $h$ lifts to
  characteristic zero.
\end{lem}

\proof The hypothesis implies that $(Z,\omega)$ is a Hurwitz tree of
type $(C,\chi)$ and conductor $h$. Hence the lemma follows from
Corollary \ref{Htreecor}. \Endproof

We set $z:=1/w$ and write
\[
  \omega=-\sum_{i=1}^\alpha \frac{2a_iz_i}{1-z_i^2w^2}{\rm d}w
        = -2\,\sum_{j\geq 0} d_jw^{2j}{\rm d}w,
\]
with
\begin{equation} \label{djeq}
     d_j = \sum_{i=1}^\alpha a_iz_i^{2j+1}.
\end{equation}
The condition that $\omega$ has a single zero at $z=\infty$ is
equivalent to the $\alpha-1$ homogenous equations
$d_0=d_1=\ldots=d_{\alpha-2}=0$. We denote by 
\[
 {\mathcal S}:= 
       \Proj(k[z_1,\ldots,z_\alpha]/(d_0, \ldots, d_{\alpha-2}))
\]
the subscheme of $\PP^{\alpha-1} $ defined by these equations.  A
point $[z_1: \cdots : z_\alpha]\in\mathcal{S}(k)$ which lies on the
hypersurface defined by the equation
\[
       \prod_i z_i \cdot \prod_{i<j}(z_i^2-z_j^2) = 0
\]
is called a {\sl trivial solution}. If $[z_1: \cdots :
z_\alpha]\in\mathcal{S}(k)$ is not a trivial solution then we call it
a {\sl good solution}. By construction, there is a bijection between
the set of good solutions and the set of Hurwitz trees of type
$(C,\chi)$, conductor $h$, discriminant $0$ and type $\mathbf{a}$.

Let $[z_1: \cdots: z_\alpha]$ be a trivial solution. Put $J=\{1, 2,
\ldots, \alpha\}$ and $J_0=\{i\in J\,|\, z_i=0\}$. Let $I\subset
J-J_0$ be a maximal subset with the property that $z_i^2\neq z_j^2$
for all $i\neq j$ in $I$. For $l\in I$,
we define $J_l=\{ i\in J\,|\, z_i^2=z_l^2\}$. Then $J$ is the disjoint
union of the subsets $J_l$.

\begin{prop}\label{trivialsolprop}
\begin{itemize}
\item[(a)] Let $[z_1: \cdots: z_\alpha]$ be a trivial solution, and
  let $(J_l)_{l\in I}$ be the corresponding decomposition of the
  indices. Then for $l\in I$ we have
  \begin{equation} \label{sumaieq}
     \sum_{i\in J_l} \nu_i a_i\equiv 0\bmod{p},
  \end{equation}
  where $\nu_i\in\{\pm 1\}$ is chosen such that $z_i=\nu_iz_l$.
\item[(b)] Suppose that we can write $J=J_0\coprod J_1\coprod \cdots
  \coprod J_n$ such that \eqref{sumaieq} holds for $l=1,\ldots,n$ and
  some choice of $\nu_i\in\{\pm 1\}$.  Then we may define a linear
  subspace ${\mathcal S}(J_l;\nu_i)\subset {\mathcal S}$ of dimension
  $n-1$, contained in the locus of trivial solutions, by putting
  $z_i:=0$ for $i\in J_0$ and $z_i:=\nu_it_l$ for $i\in J_l$ and $l=1,
  \ldots, n$.
\item[(c)] The subspace of all trivial solutions is the
  union of the subspaces ${\mathcal S}(J_l;\nu_i)$.
\end{itemize}
\end{prop}

\proof Let $[z_1: \cdots: z_\alpha]$ be a trivial solution, and let
$J_0, \ldots, J_n$ be the corresponding decomposition of the set of
indices. It is no restriction to assume that $l\in J_l$ for
$l=1,\ldots,n$. By definition, we have $d_0=\ldots=d_{\alpha-2}=0$, where
$d_j$ is defined by \eqref{djeq}. Set $d_j^{(0)}:=d_j$ and define
inductively
\[
    d_j^{(m)} := d_j^{(m-1)}-z_m^2\cdot d_{j-1}^{(m-1)},
\]
for $m=1,\ldots,n$ and $j=m,\ldots,\alpha-1$.
By induction and a straightforward computation one shows that
\[
   d_j^{(m)} = \sum_i a_i\,\big(\prod_{l=1}^m(z_i^2-z_l^2)\big)\,
                 z_i^{2(j-m)+1}.
\]
Note that the $i$th term in the above sum is zero for $i\in
J_0\cup\ldots\cup J_m$. In particular, for $m=n-1$ we get
\[
   d_j^{(n-1)} = \big(\sum_{i\in J_n} \nu_ia_i\big)\,
      (\prod_{l=1}^{n-1}(z_n^2-z_l^2)\big)\,z_n^{2(j-n)+3}.
\]
Since $z_n\neq 0$ and $z_n\neq z_l$ for $l<n$ we conclude that 
\[
     \sum_{i\in J_n} \nu_ia_i \equiv 0 \pmod{p}.
\]
For symmetry reasons, the same holds for $J_1,\ldots,J_{n-1}$. This
proves (a).

Part (b) follows by direct verification. Part (c) follows from (a) and
(b).  \Endproof

\begin{prop}\label{goodmultprop}
  Let ${\mathbf z}=[z_1: \cdots: z_\alpha]$ be a good solution.
  Then ${\mathbf z}$ is an isolated point of multiplicity one in the
  total space ${\mathcal S}$ of solutions.
\end{prop}

\proof 
Let $[z_1+\epsilon t_1: \cdots: z_\alpha+\epsilon t_\alpha]$ be
a point of ${\mathcal S}$ over the ring $k[\epsilon]$, with
$\epsilon^2=0$. Then 
\[
   \sum_i a_i(z_i+\epsilon t_i)^{2j+1} = \sum_i a_i z_i^{2j+1}
    + \epsilon\cdot(2j+1)(\sum_i a_it_i z_i^{2j}) = 0,
\]
for $j=0,\ldots,\alpha-2$. Since $2j+1<p$ and ${\mathbf z}$ is a good
solution, we conclude that
\begin{equation} \label{tieq}
     \sum_{i=1}^\alpha a_it_i z_i^{2j} = 0,
\end{equation} 
for $j=0,\ldots,\alpha-2$. By assumption we have $z_i^2\neq z_j^2$ for
$i\neq j$. Hence the Vandermonde determinant $\det(z_i^{2j})$ is
nonzero. We conclude that \eqref{tieq}, regarded as a linear equation
in the $t_i$, has a solution space of dimension $1$. Since
$t_i=z_i$ is  a solution to \eqref{tieq}, this implies
$[z_1+\epsilon t_1: \cdots: z_\alpha+\epsilon
t_\alpha]=\mathbf{z}$. The proposition follows. \Endproof

\subsection{}\label{largepsec}
To solve the local lifting problem for $h<p$ one is free to choose the
type $\mathbf{a}=(a_j)$ (\S \ref{existencesec}).  In \S
\ref{mainthmsec}, we solve the problem by showing that one particular
type works. However, it turns out that for a given $p$ most types
work. We illustrate this in this section, by using a different type to
show that the local lifting problem holds for $p$ large. As this proof
is much shorter, it may also serve as an introduction to the type of
problems one has to solve.  Let $G$ be the dihedral group of order
$2p$.

\begin{thm}\label{largepthm}
Let $\phi$ be a local $G$-action. Let $h=2\alpha+1$ be its conductor
and suppose that $p\geq 2^\alpha$. Then $\phi$ lifts to characteristic
zero.
\end{thm}

\proof
Define ${\mathbf a}=(a_1, \ldots, a_\alpha)$ inductively by 
\[
a_1=1, \qquad a_{j+1}=1+\sum_{i=1}^ja_i.
\]
Then $\sum_{i=1}^\alpha a_i=2^\alpha-1$. Therefore Proposition
\ref{trivialsolprop} together with the assumption $p\geq 2^\alpha$
implies that we do not have any trivial solutions in this case. By
Bezout's theorem, the set of good solutions is therefore nonempty.
\Endproof

\begin{rem} 
  \begin{enumerate}
  \item In fact, since the homogenous degree of $d_j$ is $2j+1$,
    Proposition \ref{goodmultprop} together with Bezout's theorem
    shows more precisely that the number of good solutions in Theorem
    \ref{largepthm} is $1\cdot 3\cdots (2\alpha-3)$.
  \item
    One can show something more general then (i).
    Suppose that, for some type ${\mathbf a}$, the locus of trivial
    solutions has dimension zero. (This can easily be read off from
    the type.) We can prove a formula for the multiplicities of the
    trivial solutions.  This gives a formula for the number of good
    solutions, which is recursive in $h$. But it is difficult to prove
    that this number is positive, although computer experiments
    suggest that this is almost always the case.
  \end{enumerate}
\end{rem}

\subsection{The proof for general $p$}\label{mainthmsec}
In this section, we prove the lifting result for small conductor and
arbitrary $p$.

\begin{thm}\label{smallhthm}
Let $\phi$ be a local $G$-action with conductor $h$. Suppose that
$h=2\alpha+1<p$. Then $f$ lifts to characteristic zero.
\end{thm}

In combination with Theorem \ref{largehthm} we obtain:

\begin{cor}\label{maincor}
The local lifting problem holds for $G=D_p$.
\end{cor}

The rest of this section concerns the proof of Theorem \ref{smallhthm}.
Put $r=\alpha-1$ and define ${\mathbf a}=(a_1=1, \ldots, a_r=1,
a_\alpha=\alpha)$. The advantage of considering this type is that
there is a large symmetry group acting: our equations are invariant
under the symmetric group acting on $z_1, \ldots, z_r$. Therefore we
may simplify the equations for the $z_i$ by considering the elementary
symmetric functions in the $z_i$. We may also suppose that
$z_\alpha=1$. The following notation replaces the previous one.

Write $h:=c_0z^r+c_1z^{r-1}+\cdots +c_r$. We consider the $c_i$ as
variables. Set $g(z):=h(-z)$ and 
\[
  f:=\frac{g}{h}\cdot\left(\frac{z-1}{z+1}\right)^\alpha, 
     \qquad
   \omega:=\frac{{\rm d} f}{f}=\frac{(d_0z^{2r}+d_1z^{2(r-1)}+\cdots
       +d_r)\,{\rm d}z}{gh(z^2-1)}.
\]
(Note that the definition of $g$ implies that $\omega$ is odd, and
hence that the numerator of $\omega$ contains only even powers of $z$.)
A straightforward computation shows that
\begin{eqnarray*} \label{delleq}
d_\ell&=&\sum_{i=0}^\ell(-1)^i c_i c_{2\ell+1-i}(2\ell+1-2i)-
  \sum_{i=0}^{\ell-1}
(-1)^ic_ic_{2\ell-1-i}(2\ell-1-2i) \\
  &+& 2\alpha\sum_{i=0}^{\ell-1} (-1)^i c_ic_{2\ell-i}
  +\alpha(-1)^\ell c_\ell^2.
\end{eqnarray*}
In particular, the $d_i$ are homogeneous of degree $2$ in $c_0,
\ldots, c_r$. We set 
\[
 {\mathcal S}:= 
       \Proj(k[c_0,\ldots,c_r]/(d_0, \ldots, d_{r-1})).
\]
A point $[c_0: \cdots : c_r]\in\mathcal{S}(k)$ which lies on the
hypersurface 
\[
      {\rm discr}(z(z^2-1)hg) = 0
\]
is  a {\sl trivial solution}. Otherwise, it is a {\sl good
  solution}. As in \S \ref{trivsolsec}, we get a bijection between
the set of good solutions and the set of Hurwitz trees of type
$(C,\chi)$, conductor $h$, discriminant $0$ and type
$\mathbf{a}=(1,\ldots,1,\alpha)$. (In contrast to \S \ref{Htree},
  we consider the marked points to be unordered in this section.)

To prove Theorem \ref{smallhthm}, it suffices to find a good solution.
In the rest of this paper, we will show that there exists in fact a
{\sl unique} good solution $[c_0:\cdots: c_r]$. This means that the
number of good solutions $[z_1:\cdots: z_\alpha]$ is  $r!=\#S_r$.
We start by describing the trivial solutions.

\begin{lem} \label{trivsollem}
  The trivial solutions are contained in the hyperplane
  $T\subset\PP^{\alpha-1}$ defined by the equation $c_0=0$. The
  hyperplane $T$ does not contain any good solutions.
\end{lem}

\proof 
If $c_0=0$ then $\deg(h)<r$. This implies immediately that $T$ does
not contain any good solutions.

Suppose that $[c_0:\cdots:c_{\alpha-1}]$ is a trivial solution
with $c_0\neq 0$. Then we may suppose that $c_0=1$. Write
$h=\prod_{i=1}^r(z-z_i)$. Then $\mathbf{z}=[z_1:\cdots:z_r:1]$ is a
trivial solution of the equations considered in \S \ref{trivsolsec}.
It follows from Proposition \ref{trivialsolprop} that there exists a
subset $J_\alpha\subset\{1,\ldots,\alpha\}$ with $\alpha\in J_\alpha$
such that
\[
    \sum_{i\in J_\alpha} \nu_i a_i = 
       \nu_\alpha\alpha + \sum_{i\in J_{\alpha}-\{\alpha\}} \nu_i 
         \equiv 0 \pmod{p}.
\]
But this is impossible, because $2\alpha<p$. By contradiction, the
lemma is proved.
\Endproof

For $\ell=0, \ldots, r-1$, let 
\[
{\mathcal S}_\ell=\{d_0=\ldots=d_\ell=0\}\subset
\PP^r=\{[c_0:c_1:\cdots:c_r]\}.
\]
Alternatively,  ${\mathcal S}_{\ell}$ corresponds to logarithmic
differentials of the above sort with a zero of order at least
${2(\ell+1)}$ at $z=\infty$. The following claim implies Theorem
\ref{smallhthm}.

\bigskip\noindent {\bf Claim:} For $\ell=0, \ldots, r-1$, let
${\mathcal S}_\ell'$ be the closure of ${\mathcal S}_\ell-T$. Then
${\mathcal S}_\ell'\subset \PP^r$ is a linear subspace of dimension
$r-\ell-1$.  In particular, $\mathcal{S}_{r-1}'$ contains a unique
good solution.

We prove the claim by induction on $\ell$. By \eqref{delleq} we have
$d_0=c_0(c_1+\alpha c_0)$. Define ${\mathcal S}_0'$ by $c_0+\alpha
c_1=0$. Then the claim holds for $\ell=0$.

Suppose that the claim holds for $\ell$. We may choose linear
parameters $t_1, \ldots, t_{r-\ell}$ for ${\mathcal S}_\ell'$ (i.e.\ 
we write $c_j$ as a linear form in the $t_i$). Since ${\mathcal
  S}_\ell'$ is not contained in the hyperplane $T$, it is no
restriction to suppose that $c_0=t_1=:t$. To prove the claim for
$\ell+1$, it suffices to show that
\begin{equation}\label{delleq2} 
   d_{\ell+1}=t(at+b), \quad\mbox{for some }
      a,b\in k[t_2, \ldots, t_{r-\ell}]\quad\mbox{  with }b\neq 0.
\end{equation}
Namely, if this holds, ${\mathcal S}_{\ell+1}'$ is obtained by
intersecting ${\mathcal S}_\ell'$ with $at+b=0$. Since $b\neq 0$, this
intersection is not completely contained in $T$. Since $d_{\ell+1}$ is
homogeneous of degree $2$, it follows that $at+b$ is linear in $t_1,
\ldots, t_{r-\ell}$. Therefore ${\mathcal S}_{\ell+1}'$ is a linear
subspace of dimension $r-\ell-2$. 

Write $\Ord_t$ for the standard valuation in $t$. We have to show that 
$\ord_t(d_{\ell+1})=1$. Suppose that $\ord_t(d_{\ell+1})\neq1$.
To deduce a contradiction, we distinguish the following cases:
\begin{itemize}
\item[(a)] $\Ord_t(d_{\ell+1})=0$,
\item[(b)] $\Ord_t(d_{\ell+1})\geq 2$.
\end{itemize}

\bigskip\noindent
{\bf Case (a)}:
Suppose that $\Ord_t(d_{\ell+1})=0$.  Write
$d_j=d_j^0+td_j^1+t^2d_j^2$ (resp.\ $c_j=c_j^0+tc_j^1$) for the
Taylor expansion of $d_j$ (resp.\ $c_j$) with respect to $t$.  It
follows from the proof of the claim for $\ell=0$ that $c_1=-\alpha
c_0=-\alpha t$, therefore $c_0^0=c_1^0=0$. Let $j$ be the minimal
$j\geq 2$ such that $c_j^0\neq 0$. The Newton polygon of $h$ shows
that there are exactly $j$ zeroes of $h$ with negative $t$-valuation
$-1/j$. Set $z=y t^{-1/j}$ and write $\omega$ as a series in $t$. Let
$\omega_1$ be the first nonzero coefficient. A short computation gives
\[
  \omega_1 = \frac{d_{\ell+1}^0\, {\rm d}
     y}{y^{2(\ell-j+1)}(y^j+c_j^0)(y^j+(-1)^j c_j^0)}.
\]
Using the fact that $\omega$ is logarithmic, it is easy to see that
$\omega_1$ is either logarithmic or exact.  Suppose that $j$ is
{even}. (The case $j$ odd is similar and left to the reader.)  Since
$\omega_1$ has double poles in the points with $y^j+c_j^0=0$, it is
not logarithmic. We claim that $\omega_1$ is not exact either. This
gives a contradiction.

After multiplying $\omega_1$ with a constant, and changing the
parameter $y$, we may write
\[
\omega_1=\frac{{\rm d}
  y}{y^{2(\ell-j+1)}(y^j-1)^2}=\frac{y^{p-2(\ell-j+1)}(y^j-1)^{p-2}\,
  {\rm d}y}{y^p(y^j-1)^p}.
\]
Put $(y^j-1)^{p-2}=\sum_{s=0}^{p-2} \delta_s y^{js}$.

 We claim that there is a unique $0\leq s\leq p-2$ with
$ p-2\ell+2j-2+js\equiv
-1\bmod{p}.$
Namely, suppose that there exist $(s,u)$ such that 
\begin{equation}\label{steq}
p-2\ell+2j-2+js=up-1.
\end{equation}
 Then $u\equiv (p-2\ell-1)/p\bmod{j}$. Since $2j\leq 2\ell\leq 2r< p$,
it follows that $1\leq u\leq j$. Therefore $u$ is uniquely determined
by (\ref{steq}). Now define 
\[
s=\frac{up+1-p+2\ell-2j}{j}\leq\frac{jp+1-p+2\ell-2j}{j}\leq p-2.
\]
This shows that there is a unique $s$ with the required properties. Note
that for this $s$ we have $\delta_s\neq 0$.

 It follows now from standard properties of the Cartier operator that
\[
\C\omega_1=\frac{\delta_s^{1/p}y^{-1+( p-2\ell+2j-2+js)/p}\, {\rm
    d}y}{y^j-1}\neq 0.
\]
This shows that $\omega_1$ is not exact. We conclude that Case (a)
does not occur.

\bigskip\noindent {\bf Case (b)}: we suppose that $\Ord_t(d_\ell)=2$.
(The case $d_\ell=0$ is similar, and left to the reader.) Then
$d_{\ell+1}=at^2$, with $a\in k^\times$. 
  
Let $j$ be as in Case (a). By assumption, the numerator of $\omega$ is
$F\, {\rm d}z$ with $F:=d_{\ell+1} z^{2(r-\ell-1)}+\cdots +d_r.$ The
proof showing that Case (a) does not occur implies that the Newton
polygon of $F$ contains a piece with slope $1/j$ which ends in the
point $(2(r-\ell-1), 2)$. Here we use that we are in Case (b).  It
follows that
\[
      d_{\ell+1}^0=\cdots=d_{\ell+j}^0=0.
\]
Let $\omega_0$ be the differential obtained from $\omega$ by
specializing $t$ to $0$. Note that $\omega_0$ corresponds to a generic
point $[0:\cdots:c_j^0:\cdots:c_r^0]$ on ${\mathcal S}'_\ell\cap T$
and that $\dim\mathcal{S}_\ell'\cap T=r-\ell-2$.
On the other hand, we have
\[
\omega_0=\frac{d_{\ell+j+1}^0z^{2(r-\ell-j-1)}+\cdots+ d_{r}^0}
       {h_0 g_0 (z^2-1)}\,{\rm d} z,
\]
where $h_0(z)=c_j^0z^{r-j}+\cdots+c_r^0$ and $g_0(z)=h_0(-z)$. Since
$c_j^0\not=0$, $\omega_0$ has a zero of order at least $2(\ell+1)$ at
$z=\infty$. Therefore, the point
$[c_j^0:\cdots:c_r^0]\subset\PP^{\tilde{r}-1}$ (with $\tilde{r}:=r-j$)
lies on the subspace $\tilde{\mathcal S}_\ell\subset\PP^{\tilde{r}-1}$
analogous to $\mathcal{S}_\ell\subset\PP^{r-1}$. Moreover, this point
does not lie on the hyperplane $\tilde{T}=\{\tilde{c}_0=0\}$. By
induction on $r$, we may assume that we have already proved our claim
for $\tilde{r}$. In particular, $\dim\tilde{\mathcal
  S}_\ell'=r-\ell-1-j<r-\ell-2$. This contradicts the assertion made
above that $\omega_0$ corresponds to a generic point of dimension
$r-\ell-2$. It follows that Case (b) does not occur either.

We conclude that $\Ord_t(\omega)=1$. Theorem \ref{smallhthm} follows.
\Endproof

\vspace{0.5cm}

\end{document}